\documentclass[11pt]{amsart}

\usepackage{geometry}
\usepackage{todonotes}
\usepackage{amsfonts, amssymb, amscd}
\numberwithin{equation}{section}

\usepackage[symbol]{footmisc}

\usepackage{bm}
\usepackage{verbatim}
\usepackage{mathrsfs}
\usepackage{graphicx}
\usepackage{tikz-cd}
\usepackage{subcaption}
\usepackage{listings}
\usepackage{subfiles}
\usepackage[toc,page]{appendix}
\usepackage{mathtools}
\usepackage{comment}
\usepackage{enumerate}
\usepackage{enumitem}
\usepackage[all]{xy}

\usepackage{graphicx}
\graphicspath{{images/}}

\usepackage{appendix}
\usepackage{hyperref}
\hypersetup{
    colorlinks=true,
    citecolor=red,
    linkcolor=blue,
    filecolor=magenta,      
    urlcolor=red,
}
\newcommand{\bb}{\bm{b}}

\newcommand{\Spec}{\mathrm{Spec}}
\newcommand{\Mm}{{\bf{M}}}

\newcommand{\fR}{\mathfrak{R}}
\newcommand{\bM}{{\bf M}}

\newcommand{\Pp}{\mathbb{P}}
\newcommand{\Qq}{\mathbb{Q}}

\newcommand{\Rr}{\mathbb{R}}

\newcommand{\Zz}{\mathbb{Z}}

\newcommand{\Span}{\operatorname{Span}}

\newcommand{\Center}{\operatorname{center}}

\newcommand{\lct}{\operatorname{lct}}

\newcommand{\Supp}{\operatorname{Supp}}

\newcommand{\mult}{\operatorname{mult}}

\newcommand{\lf}{\lfloor}
\newcommand{\rf}{\rfloor}

\newcommand{\Oo}{\mathcal{O}}
\newcommand{\Ii}{\Gamma}

\newtheorem{thm}{Theorem}[section]
\newtheorem{conj}[thm]{Conjecture}
\newtheorem{cor}[thm]{Corollary}
\newtheorem{lem}[thm]{Lemma}
\newtheorem{prop}[thm]{Proposition}

\newtheorem{claim}[thm]{Claim}

\theoremstyle{definition}
\newtheorem{defn}[thm]{Definition}

\theoremstyle{definition}
\newtheorem{rem}[thm]{Remark}

\theoremstyle{definition}

\begin{document}

\title{On effective log Iitaka fibrations and existence of complements}

\author{Guodu Chen, Jingjun Han, and Jihao Liu}

\address{Institute for Theoretical Sciences, Westlake University, Hangzhou, Zhejiang, 310024, China}
\email{chenguodu@westlake.edu.cn}

\address{Shanghai Center for Mathematical Sciences, Fudan University, Shanghai, 200438, China}
\email{hanjingjun@fudan.edu.cn}

\address{Department of Mathematics, Northwestern University, 2033 Sheridan Rd, Evanston, IL 60208, USA}
\email{jliu@northwestern.edu}

\date{\today}

\begin{abstract}
We study the relationship between Iitaka fibrations and the conjecture on the existence of complements, assuming the good minimal model conjecture. In one direction, we show that the conjecture on the existence of complements implies the effective log Iitaka fibration conjecture. As a consequence, the effective log Iitaka fibration conjecture holds in dimension $3$. In the other direction, for any Calabi-Yau type variety $X$ such that $-K_X$ is nef, we show that $X$ has an $n$-complement for some universal constant $n$ depending only on the dimension of $X$ and two natural invariants of a general fiber of an Iitaka fibration of $-K_X$. We also formulate the decomposable Iitaka fibration conjecture, a variation of the effective log Iitaka fibration conjecture which is closely related to the structure of ample models of pairs with non-rational coefficients, and study its relationship with the forestated conjectures.
\end{abstract}

\maketitle

\tableofcontents

\section{Introduction}
We work over the field of complex numbers $\mathbb C$.

Let $X$ be a smooth projective variety with non-negative Kodaira dimension. By a well-known construction of Iitaka, there exists a birational morphism $X_\infty\to X$ from a smooth projective variety $X_\infty$, and a contraction $f_\infty:X_\infty\to Z_\infty$ onto a projective variety $Z_\infty$, such that a very general fiber of $f_\infty$ is smooth with Kodaira dimension zero, and $\dim Z_\infty=\kappa(X,K_X)$. The morphism $f_\infty:X_\infty\to Z_\infty$ is referred to as an \emph{Iitaka fibration} of $K_X$ (see Definition \ref{defn: r cartier r dvisor iitaka fibration}). It is conjectured that the pluricanonical system $|mK_X|$ defines a map which is birational to an Iitaka fibration whenever the positive integer $m$ is divisible by a positive integer depending only on the dimension of $X$:

\begin{conj}[Effective Iitaka fibration, {cf. \cite[Conjecture 1.7]{HM06}})]\label{conj: ei}
Let $d$ be a positive integer. Then there exists a positive integer $m_d$ depending only on $d$, such that for any smooth projective variety $X$ of dimension $d$ with non-negative Kodaira dimension, $|mK_X|$ defines an Iitaka fibration for any positive integer $m$ divisible by $m_d$.
\end{conj}

Conjecture \ref{conj: ei} was proved when $K_X$ is big \cite{HM06,Tak06} (see also \cite{Tsu99}), when $d=2$ due to Enriques (see also \cite{Iit70}), and when $d=3$ \cite{Mor85,Kaw86,FM00,VZ07}. An important progress towards Conjecture \ref{conj: ei} is made in \cite{BZ16}, where the authors showed that there exists a positive integer $m$ depending only on $d$ and two natural invariants of the very general fibers of an Iitaka fibration of $K_X$ (the non-vanishing order and the middle Betti number), such that $|mK_X|$ defines an Iitaka fibration. Unfortunately, we don't know the boundedness of the middle Betti numbers in dimension $\geq 3$, which leaves Conjecture \ref{conj: ei} open in dimension $\geq 4$.

In practice, it is also natural to consider the following generalized version of Conjecture \ref{conj: ei}, which is known as the \emph{effective log Iitaka fibration conjecture}.

\begin{conj}[Effective log Iitaka fibration, {cf. \cite[Conjecture 1.1]{TX09}}]\label{conj: ei intro}
Let $d$ be a positive integer and $\Ii\subset[0,1]$ a DCC set. Then there exists a positive integer $m$ depending only on $d$ and $\Ii$ satisfying the following.

Assume that $(X,B)$ is an lc pair of dimension $d$ such that $B\in\Ii$ and $\kappa(X,K_X+B)\geq 0$. Then $|\lf m(K_X+B)\rf|$ defines a map which is birational to an Iitaka fibration of $K_X+B$.
\end{conj}

Conjecture \ref{conj: ei intro} was proved when $K_X+B$ is big \cite[Theorem 1.3]{HMX14}. When $\Ii\subset[0,1]\cap\Qq$ and $(X,B)$ is klt, Conjecture \ref{conj: ei intro} was proved when $d\leq 2$ \cite{Tod10}, when $d\leq 3$ and $\kappa(X,K_X)>0$ \cite{TX09,Tod10}, and when $d=4$ and $\kappa(X,K_X)=2$ \cite{TX09}.

Recently, the theory of \emph{complements}, which was introduced by Shokurov in the study of the existence of flips for threefolds \cite{Sho92}, has gradually become one of the major topics in birational geometry. This theory has played an important role in the proof of the BAB conjecture \cite{Bir19,Bir21}, the proof of the singular Yau-Tian-Donaldson conjecture and the stable degeneration conjecture \cite{Xu20,BLX22,LXZ22,XZ22}, and recent studies on Shokurov's ascending chain condition conjecture for minimal log discrepancies \cite{Liu18,HLS19,HL22b,HLL22}. For other related work, we refer the readers to \cite{Sho20,CGN21,CZ21,CHX22,CX22,Liu22}. Although the existence of $n$-complements \cite{Bir19,HLS19} is settled for Fano type varieties, it is natural to consider the existence of $n$-complements for pairs admitting an lc Calabi-Yau structure, that is, \emph{$\Rr$-complementary varieties}. We have the following conjecture:

\begin{conj}[Existence of complements]\label{conj: boundedness and existence of n complement nft}
Let $d$ be a positive integer and $\Ii\subset [0,1]$ a DCC set. Then there exists a positive integer $n$ depending only on $d$ and $\Ii$ satisfying the following. 

Assume that $(X/Z\ni z,B)$ is an $\Rr$-complementary pair of dimension $d$ such that $B\in\Ii$, then $(X/Z\ni z,B)$ has an $n$-complement $(X/Z\ni z,B^+)$. Moreover, if the closure of $\Ii$ belongs to $[0,1]\cap\mathbb Q$, then we can pick $B^+\ge B$.
\end{conj}

\medskip

\noindent\textbf{Relationship between the conjectures}. It is interesting to ask whether we can establish some connections between Conjectures \ref{conj: ei}, \ref{conj: ei intro} and Conjecture \ref{conj: boundedness and existence of n complement nft}. At first glance, it is difficult to observe their relationships, as these conjectures are considering structures of varieties and pairs with completely different positivity properties: Conjectures \ref{conj: ei}, \ref{conj: ei intro} are concentrating on pairs with positive canonical bundle (i.e., $K_X+B$ is effective), while Conjecture \ref{conj: boundedness and existence of n complement nft} is concentrating on varieties with negative canonical bundle (i.e., $-(K_X+B)$ is effective). Surprisingly, we have the following theorems, which show that these conjectures are actually deeply related with each other in multiple directions.

\medskip

\noindent\textit{From complements to Iitaka fibrations}. First, we prove that Conjecture \ref{conj: boundedness and existence of n complement nft} implies Conjecture \ref{conj: ei intro} assuming the good minimal model conjecture.

\begin{thm}\label{thm: main}
Let $d$ be a positive integer. Assume that the good minimal model conjecture (Conjecture \ref{conj: exist gmm}) and the existence of complements (Conjecture \ref{conj: boundedness and existence of n complement nft}) hold in dimension $d$. Then the effective log Iitaka fibration conjecture (Conjecture \ref{conj: ei intro}) holds in dimension $d$. 
\end{thm}

As an immediate corollary, we prove Conjecture \ref{conj: ei intro} in dimension $\leq 3$:

\begin{cor}\label{cor: ei dim 3}
 Conjecture \ref{conj: ei intro} holds when $d\leq 3$. 
\end{cor}

\noindent\textit{From Iitaka fibrations to complements}. Next, we show that the existence of $n$-complements is deeply related to some invariants associated to Iitaka fibrations. We have the following theorem:

\begin{thm}\label{thm: betti to complement intro}
Let $d,b,$ and $\beta$ be three positive integers. Assume that the good minimal model conjecture (Conjecture \ref{conj: exist gmm}) holds in dimension $d$. Then there exists a positive integer $n$ depending only on $d,b,$ and $\beta$ satisfying the following. 

Assume that $X$ is a $\Qq$-factorial normal projective variety of dimension $d$, $f_\infty:X_\infty\to Z_\infty$ is an Iitaka fibration of $-K_X$, $h:X_\infty\to X$ is the induced birational morphism, and $F$ is a general fiber of $f_\infty$. Suppose that
\begin{enumerate}
    \item $-K_X$ is nef,
    \item $X$ has a klt $\Rr$-complement,
    \item $\kappa(X,-K_X)\ge0$, and $b$ is the non-vanishing order of $-h^*K_X|_F$, i.e.,
    $$b=\min\left\{a\in\Zz_{>0}\mid |-ah^*K_X|_F|\not=\emptyset\right\},$$
and
    \item $\beta:=\dim H^{\dim F}(\tilde F,\mathbb C)$, where $\tilde F$ is a smooth model of the cover of $F$ associated to the unique divisor of $|-bh^*K_X|_F|$. 
\end{enumerate}
Then $X$ has an $n$-complement.
\end{thm}
We remark that the assumptions on $b$ and $\beta$  in Theorem \ref{thm: betti to complement intro} are natural assumptions, and are exactly the additional assumptions in \cite[Theorem 1.2]{BZ16} on the effective Iitaka fibration conjecture. We also note that, modulo the good minimal model conjecture, the boundedness of $b$ follows immediately from the effective log Iitaka conjecture (Conjecture \ref{conj: ei intro}) for log pairs with Iitaka dimension $0$.

\medskip

\noindent\textbf{Decomposable Iitaka fibrations}. For any semi-ample divisor $D$, the ample model of $D$ is clearly birational to an Iitaka fibration of $D$. However, if $D$ is only assumed to be an $\Rr$-divisor, then it is possible that $D$ does not have any Iitaka fibration although the ample model of $D$ exists \cite[Example 1.2]{CHL22}. An approach to resolve this issue is to define the \emph{invariant Iitaka fibration} (see Definition \ref{defn: invariant iitaka fibration}). A question arises naturally: \emph{do we expect any kind of uniform effectively on invariant Iitaka fibrations for log pairs}, similar to the effective log Iitaka fibration conjecture? 

The question is expected to have a positive answer. More precisely, suppose that $(X,B)$ is an lc pair such that $K_X+B$ induces a map $X\dashrightarrow Z$ that is birational to an invariant Iitaka fibration of $K_X+B$. Then we expect that $X\dashrightarrow Z$ will actually be birational to an Iitaka fibration of $K_X+B'$ after we (uniformly) perturb the coefficients of the boundary $B$ to get a new boundary $B'$. Therefore, the effective log Iitaka fibration conjecture should induce some kind of uniform effectivity on an invariant Iitaka fibration induced by $K_X+B$.

The difficulty is to show that we can make a uniform perturbation to switch the boundary $B$ to a new boundary $B'$. In \cite[Theorem 1.1]{CHL22}, the authors prove a weaker version, which shows that, modulo the non-vanishing conjecture, there exists such uniform perturbation such that an Iitaka dimension of $K_X+B'$ is equal to an invariant Iitaka dimension of $K_X+B$. In this paper, we will show that an Iitaka fibration of $K_X+B'$ is actually equal to an Iitaka fibration as $K_X+B$. To make our statements more clear, we introduce the concept of \emph{decomposable Iitaka fibrations}.

\begin{defn}[Decomposable Iitaka fibrations]
Let $\Ii_0:=\{a_1,\dots,a_k\}\subset (0,1]$ be a finite set such that $\sum_{i=1}^ka_i=1$, and $\Ii'\subset[0,1]$ a set. Assume that $(X,B)$ is an lc pair such that $\kappa_\iota(X,K_X+B)\ge0$. We say that $(X,B)$ has a \emph{$(\Ii_0,\Ii')$-decomposable Iitaka fibration} if there exist $\Rr$-divisors $B_1,\dots,B_k\in\Ii'$ such that
\begin{enumerate}
    \item $B=\sum_{i=1}^k a_iB_i$,
    \item $(X,B_i)$ is lc for any $1\le i\le k$, and
    \item any invariant Iitaka fibration of $K_X+B$ is an Iitaka fibration of $K_X+B_i$ for any $1\le i\le k$.
\end{enumerate}
In addition, if
\begin{enumerate}\setcounter{enumi}{3}
    \item the map defined by $|\lf m(K_X+B_i)\rf|$ is birational to any invariant Iitaka fibration of $K_X+B$ for any integer $1\le i\le k$,
\end{enumerate}
then we say that $(X,B)$ has an \emph{$(m,\Ii_0,\Ii)$-decomposable Iitaka fibration}. 
\end{defn}

As an analogue of the effective log Iitaka fibration conjecture, we propose the following conjecture on the existence of decomposable Iitaka fibrations. Roughly speaking, the conjecture indicates that we can uniformly perturb an invariant Iitaka fibration to get an effective log Iitaka fibration.

\begin{conj}[Decomposable  Iitaka fibrations]\label{conj: decomposable ei intro}
Let $d$ be a positive integer and $\Ii\subset[0,1]$ a DCC set. Then there exist a positive integer $m$, a finite set $\Ii_0\subset (0,1]$, and a DCC set $\Ii'\subset [0,1]$ depending only on $d$ and $\Ii$ satisfying the following. Assume that $(X,B)$ is a $\Qq$-factorial lc pair of dimension $d$ such that $B\in\Ii$ and $\kappa_\iota(X,K_X+B)\ge0$. Then: \begin{enumerate}
    \item (Weak version) $(X,B)$ has a $(\Ii_0,\Ii')$-decomposable Iitaka fibration.
    \item (Strong version) $(X,B)$ has an $(m,\Ii_0,\Ii')$-decomposable Iitaka fibration.
\end{enumerate}
\end{conj}
Conjecture \ref{conj: decomposable ei intro} will help us to understand the structure of good minimal models and their ample models for pairs with real coefficients. We show that Conjecture \ref{conj: decomposable ei intro} follows from the non-vanishing conjecture (Conjecture \ref{conj: non-vanishing}) and the effective log Iitaka fibration conjecture (Conjecture \ref{conj: decomposable ei intro}(2)).

\begin{thm}\label{thm:dif1}
Let $d$ be a positive integer. Assume that the non-vanishing conjecture (Conjecture \ref{conj: non-vanishing}) holds in dimension $d$. Then:
\begin{enumerate}
\item Conjecture \ref{conj: decomposable ei intro}(1) holds in dimension $d$.
\item Assume that Conjecture \ref{conj: ei intro} holds in dimension $d$. Then  Conjecture \ref{conj: decomposable ei intro}(2) holds in dimension $d$.
\end{enumerate}
\end{thm}

As an immediate corollary, we have:

\begin{cor}\label{cor: decomposable ei dim 3}
Conjecture \ref{conj: decomposable ei intro} holds when $d\leq 3$.
\end{cor}

Combining Theorems \ref{thm: main} and \ref{thm:dif1}, we show that Conjecture \ref{conj: decomposable ei intro} follows from the good minimal model conjecture and the existence of complements.

\begin{thm}\label{thm:dif2}
Let $d$ be a positive integer. Assume that the good minimal model conjecture (Conjecture \ref{conj: exist gmm}) and the existence of complements (Conjecture \ref{conj: boundedness and existence of n complement nft}) hold in dimension $d$. Then Conjecture \ref{conj: decomposable ei intro} holds in dimension $d$.
\end{thm}

Finally, we remark that we expect all theorems in our paper to hold in the relative setting. That is, instead of considering Iitaka fibrations for projective varieties and projective pairs, we may also consider Iitaka fibrations in the relative case (cf. \cite[Definition 3.19]{Li22}). 

\medskip

\noindent\textbf{Structure of the paper}. In Section \ref{sec2}, we introduce some notation and tools which will be used in this paper. In Section \ref{sec3}, we recall the canonical bundle formulas and prove Proposition \ref{prop:cbfindex}. In Section \ref{sec4}, we prove Theorem \ref{thm: main}. In Section \ref{sec6}, we prove Theorem \ref{thm: betti to complement intro}. In Section \ref{sec5}, we introduce invariant Iitaka fibrations and prove Theorems \ref{thm:dif1} and \ref{thm:dif2}.

\medskip

\noindent\textbf{Acknowledgement}.
The second and the third named authors began this work when they worked on \cite{HL20} in Summer 2020. Part of this work was done while the first named author visited Chuyu Zhou at EPFL in Summer 2022. He would like to thank their hospitality. The authors would like to thank Qianyu Chen, Christopher D. Hacon, Fei Hu, Chen Jiang, Junpeng Jiao, Zhan Li, Haidong Liu, Wenfei Liu, Yuchen Liu, Yujie Luo, Fanjun Meng, Lingyao Xie, Qingyuan Xue, and Chuyu Zhou for valuable discussions and suggestions. The authors would like to thank Enrica Floris for answering questions about Theorem \ref{thm: betti to complement intro}. The first named author was supported by the China post-doctoral grants BX2021269 and 2021M702925. The second named author was supported by National Key Research and
Development Program of China (Grant No. 2020YFA0713200). The second named author is a member of LMNS, Fudan University.

\section{Preliminaries}\label{sec2}
We adopt the standard notation and definitions in \cite{KM98,BCHM10} and will freely use them. 

\subsection{Divisors}

\begin{defn}\label{defn: DCC and ACC}
Let $\Ii\subset\Rr$ be a set. We say that $\Ii$ satisfies the \emph{ascending chain condition} (ACC) if any increasing sequence in $\Ii$ stabilizes. We say that $\Ii$ satisfies the \emph{descending chain condition} (DCC) if any decreasing sequence in $\Ii$ stabilizes. 
\end{defn}

\begin{defn}[{\cite[3.2]{PS09}, \cite[2.2]{Bir19}}]
Let $\fR\subset[0,1]\cap\Qq$ be a finite set, we define
$$\Phi(\fR):=\left\{1-\frac{\gamma}{n}\mid\gamma\in\fR,n\in\Zz_{\ge 1}\right\}.$$
We say that a set $\Ii\subset[0,1]$ is a \emph{hyperstandard set} if there exists a finite set $\fR\subset [0,1]\cap\Qq$ such that $0,1\in\fR$ and $\Ii=\Phi(\fR)$.
\end{defn}

\begin{defn}
We say $f: X \to Z$ is a \emph{contraction} if $\pi$ is a projective morphism, and $f_*\Oo_X = \Oo_Z$. We say that a birational map $\phi: X \dashrightarrow Y$ is a \emph{birational contraction} if $\phi$ is projective and $\phi^{-1}$  does not contract any divisors.
\end{defn}

\begin{defn}
Let $\Ii\subset\Rr$ be a set, $X$ a variety, and $B:=\sum_{i=1}^s b_iB_i$ an $\Rr$-divisor on $X$, where $B_i$ are the irreducible components of $B$. We write $B\in\Ii$ if $b_i\in\Ii$ for every $i$. We define 
$$|B|:=\max_{1\le i\le s}{|b_i|}, \lf B\rf:=\sum_{i=1}^s \lf b_i\rf B_i,\text{ and }\{ B\}:=\sum_{i=1}^s \{ b_i\} B_i.$$
We may denote by $\mathcal{K}(X)$ the rational function field of $X$.

Let $f: X\rightarrow Z$ be a contraction between normal quasi-projective varieties and $D$ an $\Rr$-divisor on $X$. We say that $D$ is \emph{vertical} over $Z$ if $f(\Supp D)$ is a proper subset of $Z$. We say that $D$ is \emph{horizontal} over $Z$ if $D$ is not vertical over $Z$. We may uniquely write $D=D^h+D^v$ such that $D^v$ is vertical over $Z$ and each component of $D^h$ is horizontal over $Z$. We call $D^h$ the \emph{horizontal$/Z$ part} of $D$ and $D^v$ the \emph{vertical$/Z$ part} of $D$.
\end{defn}

\begin{lem}\label{lem:efftri}
Let $f:X\to Z$ be a contraction between normal quasi-projective varieties and $L$ a Cartier divisor on $X$ such that $L\sim_{Z}0$. Suppose either $L$ is vertical over $Z$ or $L\ge0$. Then there is a Cartier divisor $L_Z$ on $Z$ such that $L=f^*L_Z$.
\end{lem}

\begin{proof}
If $L\ge0$, then $L$ is vertical over $Z$. Thus it suffices to prove the statement under the condition that $L$ is vertical over $Z$. By assumption there exist a rational function $s\in\mathcal{K}(X)$ and a Cartier divisor $L_Z'$ on $Z$ such that $L-f^*L_Z'=(s).$ Since $L$ is vertical over $Z$, we may find an open subset $V\subset Z$ such that $(s)|_U\ge0$, where $U:=f^{-1}V\subset X$. In particular, $s\in \Oo_X(U)$ and thus $s\in \Oo_Z(V)$ as $f$ is a contraction and $\Oo_Z(V)=(f_*\Oo_X)(V)=\Oo_X(U).$ Hence $s\in \mathcal{K}(V)=\mathcal{K}(Z)\hookrightarrow 
\mathcal{K}(X)$ and $\phi(s)=s$, where $\phi:\mathcal{K}(Z)\hookrightarrow \mathcal{K}(X)$. Thus $f^*(s)=(\phi(s))=(s)$. Set $L_Z:=L_Z'+(s)$ and we are done.
\end{proof}

\begin{lem}\label{lem:trivialoverbase}
Let $f:X\to Z$ be a contraction between normal quasi-projective varieties, $D$ a divisor on $X$, and $U_1,U_2\subset Z$ two open subsets such that $D|_{f^{-1}(U_i)}\sim_{U_i}0$ for $i=1,2$. Then $D|_{f^{-1}(U_1\cup U_2)}\sim_{U_1\cup U_2}0$.
\end{lem}
\begin{proof}
Replacing $Z$ by $U_1\cup U_2$, we may assume that $Z=U_1\cup U_2$. By assumption, 
$$D|_{f^{-1}(U_i)}= f_i^*D_i+(s_i)$$
for some Cartier divisor $D_i$ on $U_i$ and $s_i\in \mathcal{K}(X)$, where $f_i:=f|_{f^{-1}(U_i)}$ for $i=1,2$. Let $f_{12}:=f|_{f^{-1}{(U_1\cap U_2)}}$. Then $(f_{12})^*(D_1-D_2)=(\frac{s_2}{s_1})$ on $f^{-1}(U_1\cap U_2)$. By the projection formula, we have
\begin{align*}
\mathcal{O}_{U_1\cap U_2}&\cong 
(f_{12})_{*}\left((\frac{s_2}{s_1})\mathcal{O}_{X}|_{f^{-1}(U_1\cap U_2)}\right)\\&=(f_{12})_*(f_{12})^*\mathcal{O}_{U_1\cap U_2}(D_1-D_2)\\&=\mathcal{O}_{U_1\cap U_2}(D_1-D_2).
\end{align*}
Thus $(D_1-D_2)|_{U_1\cap U_2}=(s_Z)$ and $(f_{12})^*(s_Z)=(\frac{s_2}{s_1})$ over $U_1\cap U_2$ for some $s_Z\in \mathcal{K}(Z)$. In particular, $(\frac{s_2}{s_1})-f^*(s_Z)$ is vertical over $Z$. Note that $(\frac{s_2}{s_1})-f^*(s_Z)\sim_{Z}0$. By Lemma \ref{lem:efftri}, there exists a Cartier divisor $L_Z$ on $Z$ such that 
$$(\frac{s_2}{s_1})-f^*(s_Z)=f^{*}L_Z.$$ 
Moreover, we have $\Supp L_Z\cap (U_1\cap U_2)=\emptyset$, and $L_Z+(s_Z)=D_1-D_2$ on $U_1\cap U_2$ as $(D_1-D_2)|_{U_1\cap U_2}=(s_Z)$. Hence there exists a Cartier divisor $D'$ on $Z$, such that
$$D'=D_1\text{ on }U_1\text{ \ and \ }D'=D_2+(s_Z)+L_Z\text{ on }U_2.$$
It follows that $f^*D'=f^*D_1=D-(s_1)$ over $U_1$ and $f^*D'=f^*(D_2+(s_Z)+L_Z)=D-(s_2)+(\frac{s_2}{s_1})=D-(s_1)$ over $U_2$. Hence $D=f^*D'-(s_1)$, and thus $D\sim_{Z} 0$.
\end{proof}

\subsection{Pairs and singularities}

\begin{defn}[Pairs, {cf. \cite[Definition 2.2]{CH21}}] \label{defn sing}
A \emph{pair} $(X/Z\ni z, B)$ consists of a contraction $f: X\rightarrow Z$ between normal quasi-projective varieties, a (not necessarily closed) point $z\in Z$, and an $\mathbb{R}$-divisor $B\geq 0$ on $X$, such that $K_X+B$ is $\Rr$-Cartier over a neighborhood of $z$. We may also call it an \emph{$\Rr$-pair}. If $B\in\Qq$, then we call $(X/Z\ni z, B)$ a \emph{$\Qq$-pair}. If $f=id$ and $z=x\in X$, then we may use $(X\ni x, B)$ instead of $(X/Z\ni z,B)$. If $(X/Z\ni z,B)$ (resp., $(X\ni x,B)$) is a pair for any point $z\in Z$ (resp., $x\in X$), then we call $(X/Z,B)$ (resp., $(X,B)$) a pair.
\end{defn}
\begin{defn}[Singularities of pairs]\label{defn: relative mld}
 Let $(X/Z\ni z,B)$ be a pair associated with the contraction $f: X\to Z$, and let $E$ be a prime divisor over $X$ such that $z\in f(\Center_X E)$. Let $g: Y\rightarrow X$ be a log resolution of $(X,B)$ such that $\Center_Y E$ is a divisor, and suppose that $K_Y+B_Y=g^*(K_X+B)$ over a neighborhood of $z$. We define $a(E,X,B):=1-\mult_EB_Y$ to be the \emph{log discrepancy} of $E$ with respect to $(X,B)$. 
 
We say that a prime divisor $E$ is \emph{over} $X/Z\ni z$ if $E$ is a prime divisor $E$ over $X$ and $f(\Center_X E)=\bar z$. We say that $(X/Z\ni z,B)$ is \emph{lc} (resp., \emph{klt}) if $a(E,X,B)\geq 0$ (resp., $>0$) for any prime divisor $E$ over $X/Z\ni z$. We say that $(X/Z,B)$ is \emph{lc} (resp., \emph{klt}) if $(X\ni x,B)$ is lc (resp., klt) for any codimension $\geq 1$ point $x\in X$. We say that $(X/Z,B)$ is \emph{dlt} if there exists a log resolution $g: Y\rightarrow X$ of $(X,B)$ such that $a(E,X,B)>0$ for any $g$-exceptional prime divisor $E$.
\end{defn} 

\begin{lem}\label{lem:gfiso}
Assume that $(X/Z,B)$ is an lc pair and $D$ is an $\Rr$-Cartier $\Rr$-divisor on $X$ such that $K_X+B\sim_{\Rr,Z} D$ and $D$ is vertical over $Z$. Suppose that $\psi:X\dashrightarrow X'$ is a partial $(K_X+B)$-MMP over $Z.$ Then there exists a proper open subset $U\subset Z$ such that $X_z$ is isomorphic to $X'_z$ for any $z\in U$, where $X_z$ (resp., $X'_z$) is the fiber of $X\to Z$ (resp., $X'\to Z$) over $z$.
\end{lem}
\begin{proof}
Let $f: X\to Z$ be the associated morphism. Since $K_X+B\sim_{\mathbb R}0$ over $Z\backslash f(\Supp D)$, $\psi$ is an isomorphism over $Z\backslash f(\Supp D)$, so we may choose $U:=Z\backslash f(\Supp D)$.
\end{proof}

\begin{lem}\label{lem: dim 1 base q divisor}
Assume that $(X/Z\ni z,B)$ is an lc pair such that $Z$ is a curve, $K_X+B\sim_{\mathbb R}0$ over a neighborhood of $z$, $B^h\in\Qq$, and $z$ is a closed point. Then $B+sf^*z$ is a $\Qq$-divisor over a neighborhood of $z$, where $s:=\lct(X/Z\ni z,B;f^*z)$ and $f: X\rightarrow Z$.
\end{lem}
\begin{proof}
Possibly shrinking $Z$ we may assume that $K_X+B\sim_{\Rr,Z}0$ and $f(\Supp B^v)=\{z\}$. There exist real numbers $r_1,\dots,r_c$, $\Qq$-linear functions $s_1,\dots,s_m:\mathbb R^{c+1}\rightarrow\mathbb R$, and Weil divisors $B_1,\dots,B_m$ on $X$, such that $1,r_1,\dots,r_c$ are linearly independent over $\Qq$, and $B=\sum_{i=1}^ms_i(1,\bm{r}_0)B_i$, where $\bm{r}_0:=(r_1,\dots,r_c)$. Let $B(\bm{v}):=\sum_{i=1}^ms_i(1,\bm{v})B_i$ for any $\bm{v}\in\mathbb R^c$. Since $B^h\in\mathbb Q$, $B^h=B(\bm{v})^h\in\Qq$ for any $\bm{v}\in\Rr^c$.

By \cite[Lemma 5.3]{HLS19}, $K_X+B(\bm{v})\sim_{\mathbb R,Z}0$ for any $\bm{v}\in\mathbb R^c$ and thus $B(\bm{v})-B\sim_{\mathbb R,Z}0$. Moreover, as $f(\Supp(B(\bm{v})-B))=f(\Supp B^v)=\{z\}$, we see that $B-B(\bm{v})=l_{\bm{v}}f^*z$ for some real number $l_{\bm{v}}$. Pick $\bm{v}_0\in\mathbb Q^c$, then
$$s+l_{{\bm{v}}_0}=\lct(X/Z\ni z,B(\bm{v}_0);f^*z)\in\mathbb Q$$
 as $B(\bm{v}_0)\in\mathbb Q$. Since
$$B+sf^*z=B(\bm{v_0})+(s+l_{{\bm{v}}_0})f^*z,$$
$B+sf^*z$ is a $\Qq$-divisor. This finishes the proof.
\end{proof}

\begin{conj}[Good minimal model conjecture]\label{conj: exist gmm}
Let $d$ be a positive integer. Assume that $(X/Z,B)$ is an lc pair of dimension $d$ such that $K_X+B$ is pseudo-effective over $Z$. Then $(X/Z,B)$ has a good minimal model over $Z$.
\end{conj}

\begin{conj}[Non-vanishing conjecture]\label{conj: non-vanishing}
Let $d$ be a positive integer. Assume that $(X/Z,B)$ is an lc pair of dimension $d$ such that $K_X+B$ is pseudo-effective over $Z$. Then $|(K_X+B)/Z|_{\Rr}\not=\emptyset$. 
\end{conj}

\subsection{Complements}

\begin{defn}\label{defn complement}
Let $n$ be a positive integer, $\Ii\subset (0,1]$ a set, and $(X/Z\ni z,B)$ and $(X/Z\ni z,B^+)$ two pairs. We say that $(X/Z\ni z,B^+)$ is an \emph{$\Rr$-complement} of $(X/Z\ni z,B)$ if $(X,B^+)$ is lc, $B^+\geq B$, and $K_X+B^+\sim_{\Rr}0$ over a neighborhood of $z$. We say that $(X/Z\ni z,B)$ is \emph{$\Rr$-complementary} if $(X/Z\ni z,B)$ has an $\Rr$-complement. 

We say that $(X/Z\ni z,B^+)$ is an \emph{$n$-complement} of $(X/Z\ni z,B)$ if
\begin{itemize}
\item $(X/Z\ni z,B^+)$ is lc,
\item $nB^+\geq \lfloor (n+1)\{B\}\rfloor+n\lfloor B\rfloor$, and
\item $n(K_X+B^+)\sim 0$ over a neighborhood of $z$.
\end{itemize}
We say that $(X/Z\ni z,B^+)$ is an \emph{$(n,\Ii)$-decomposable $\Rr$-complement} of $(X/Z\ni z,B)$ if there exist real numbers $a_1,\dots,a_k\in\Ii$, and $\Qq$-divisors $B_1^+,\dots,B_k^+$ on $X$, such that
\begin{itemize}
\item $\sum_{i=1}^ka_i=1$ and $\sum_{i=1}^ka_iB_i^+=B^+$,
\item $(X/Z\ni z,B^+)$ is an $\Rr$-complement of $(X/Z\ni z,B)$, and
\item  $(X/Z\ni z,B_i^+)$ is an $n$-complement of itself for any integer $1\le i\le k$.
\end{itemize}
\end{defn} 

Conjecture \ref{conj: boundedness and existence of n complement nft} holds true when $d=3$. More precisely, we have the following.
\begin{thm}\label{thm:3folddcc}
Let $l$ be a positive integer and $\Ii\subset[0,1]$ a DCC set. Then there exists a positive integer $n$ which is divisible by $l$ depending only on $l$ and $\Ii$ satisfying the following.

Assume that $(X/Z\ni z,B)$ is an $\Rr$-complementary pair of dimension $3$ with $B\in\Ii$. Then $(X/Z\ni z,B)$ has an $n$-complement $(X/Z\ni z,B^+)$. Moreover, if $\Span_{\Qq_{\ge0}}(\bar{\Ii}\backslash\Qq)\cap (\Qq\backslash\{0\})=\emptyset$, then we can pick $B^+\ge B.$ 
\end{thm}
\begin{proof}
This follows from \cite[Theorem 1]{FMX19} and \cite[Theorem 8.25]{HLS19}.
\end{proof}

\subsection{Iitaka dimensions and invariant Iitaka dimensions}

\begin{defn}[Iitaka dimensions, cf. {\cite[II 3.2 Definition]{Nak04}}]\label{defn: Iitaka dimension}
Let $X$ be a normal projective variety and $D$ an $\Rr$-divisor on $X$. For any positive integer $m$ such that $|\lfloor mD\rfloor|\not=\emptyset$, we denote $\Phi_m: X\dashrightarrow\Pp(H^0(X,\lfloor mD\rfloor)).$ The \emph{Iitaka dimension} $\kappa(X,D)$ of $D$ is defined in the following way. If $|\lfloor mD\rfloor|\not=\emptyset$ for some positive integer $m$, then
$$\kappa(X,D):=\max\{\dim \Phi_m(X)\mid {m\in\Zz_{>0},|\lfloor mD\rfloor|\not=\emptyset}\}.$$
Otherwise, let $\kappa(X,D):=-\infty$. Note that if $|\lfloor mD\rfloor|\not=\emptyset$, then by \cite[II 3.8 Corollary]{Nak04},
$$\kappa(X,D)=\max\left\{k\in\Zz_{>0}\mid \underset{m\rightarrow+\infty}{\lim\sup}\frac{\dim H^0(X,\lfloor mD\rfloor)}{k^m}>0\right\}.$$
\end{defn}

\begin{defn}[Invariant Iitaka dimensions, cf. {\cite[Definition 2.2.1]{Cho08}}]
Let $X$ be a normal projective variety and $D$ an $\Rr$-divisor on $X$. The \emph{invariant Iitaka dimension} $\kappa_{\iota}(X,D)$ of $D$ is defined as follows. If $|D|_{\Rr}\not=\emptyset$, then we define
$$\kappa_{\iota}(X,D):=\kappa(X,D')$$
for some $\Rr$-divisor $D'\in|D|_\Rr$. Otherwise, we let $\kappa_{\iota}(X,D):=-\infty$. Note that $\kappa_{\iota}(X,D)$ is independent of the choice of $D'$ \cite[Corollary 2.1.4]{Cho08}.
\end{defn}

We gather some basic properties of $\kappa$ and $\kappa_{\iota}$ which will be used in the rest of this paper.

\begin{prop}[{\cite[Propositions 2.1.2, 2.2.2, Corollary 2.1.4]{Cho08}}]\label{prop: basic properties of three iitaka dimensions}
Let $X$ be a normal projective variety, and $D\sim_\Rr D'$ two $\Rr$-Cartier $\Rr$-divisors on $X$. Then
\begin{enumerate}
  \item $\kappa(X,D)\leq\kappa_{\iota}(X,D)=\kappa_{\iota}(X,D')$, and $\kappa(X,D)<\kappa_{\iota}(X,D)$ if and only if $\kappa(X,D)=-\infty$ and $\kappa_{\iota}(X,D)\geq 0$.
  \item If $D'\geq 0$, then $\kappa(X,D)\leq\kappa(X,D')$. 
\end{enumerate}
\end{prop}

\cite[Example 2.10]{CHL22} shows that we may have strict inequality in Proposition \ref{prop: basic properties of three iitaka dimensions}(1). 
\begin{lem}[{\cite[Proposition 3.20]{Sho03},{~\cite[II Lemma 3.11]{Nak04}}}]\label{lem: pullback three iitaka dimensions}
Let $g: Y\rightarrow X$ be a surjective morphism between normal projective varieties and $D$ an $\Rr$-Cartier $\Rr$-divisor on $X$. Then
\begin{enumerate}
  \item $\kappa(X,D)=\kappa(Y,g^*D)$ and $\kappa_{\iota}(X,D)=\kappa_{\iota}(Y,g^*D)$, and
  \item if $g$ is birational, then $\kappa(X,D)=\kappa(Y,g^*D+E)$ and $\kappa_{\iota}(X,D)=\kappa_{\iota}(Y,g^*D+E)$ for any $g$-exceptional $\Rr$-Cartier $\Rr$-divisor $E\geq 0$ on $Y$.
\end{enumerate}
\end{lem}

\subsection{Iitaka fibrations for $\Rr$-divisors}

\begin{defn}
Assume that $f:X\dashrightarrow Z$ is a rational map and $f_\infty: X_{\infty}\rightarrow Z_{\infty}$ is a projective morphism. We say that \emph{$f$ is birational to $f_\infty$} if there exist a birational morphism $h: X_{\infty}\rightarrow X$ and a birational map $g:Z_{\infty}\dashrightarrow Z'$ such that $f\circ h=g\circ f_\infty$, i.e., the following diagram is commutative:
    \begin{center}$\xymatrix{
    X_{\infty} \ar@{->}[d]_{f_\infty}\ar@{->}[rr]^{h}   &  & X\ar@{-->}[d]^{f} \\
     Z_{\infty}\ar@{-->}[rr]^{g}& &Z. }$
\end{center}
\end{defn}

\begin{defn}[cf. {\cite[Definition 3.19]{Li22}}]\label{defn: r cartier r dvisor iitaka fibration}
Let $X$ be a normal projective variety and $D$ an $\Rr$-Cartier $\Rr$-divisor on $X$ such that $\kappa(X,D)\geq 0$. 
A projective morphism $f_\infty: X_{\infty}\rightarrow Z_{\infty}$ between quasi-projective smooth varieties is called an \emph{Iitaka fibration} of $D$ if the following hold:
\begin{enumerate}
    \item $\dim Z_{\infty}=\kappa(X,D)$,
    \item $f_m: X\dashrightarrow Z_m\subset\mathbb PH^0(X,\lfloor mD\rfloor)$, the map associated with the complete linear system $|\lf mD\rf|$, is birational to $f_\infty$ with the morphism $h:X_\infty\to X$, for any sufficiently divisible large integer $m$, and
    
    \item for any sufficiently large integer $n$, we have 
    $$\kappa\left(F,\left(h^{-1}_*D+nE\right)|_F\right)=0,$$
    where $F$ is a very general fiber of $f_\infty$ and $E$ is the sum of all the $h$-exceptional prime divisors.
\end{enumerate}
\end{defn}

By \cite[Proposition 3.20]{Li22}, for any $\Rr$-Cartier $\Rr$-divisor $D$ such that $\kappa(X,D)\geq 0$, there always exists an Iitaka fibration of $D$. The following lemmas are well-known to experts.

\begin{lem}\label{lem:replaceif0}
Notation as in Definition \ref{defn: r cartier r dvisor iitaka fibration}.
\begin{enumerate}
  \item Assume that $D_\infty$ is an $\Rr$-divisor on $X_\infty$ such that $D_\infty-h^*D\ge0$ and is $h$-exceptional, then $f_\infty$ is an Iitaka fibration of $D_\infty.$
  \item Assume that $X_\infty'\to X_\infty$ is a birational morphism from a smooth variety, then $X_\infty'\to Z_\infty$ is an Iitaka fibration of $D$.
\end{enumerate}
\end{lem}
\begin{proof}
(1) By assumption, $H^0(X,\lf mD\rf)=H^0(X_\infty,\lf mD_\infty\rf)$ and $X_\infty\dashrightarrow Z_m$ is also the map associated with the complete linear system $|\lf mD_\infty\rf|$. Then (1) follows from the definition. 

(2) We only need to show $\kappa(F',((h')^{-1}_*D+nE')|_{F'})=0$ for any sufficiently large integer $n$, where $h':X_\infty'\to X$, $F'$ is a very general fiber of $X_\infty'\to Z_\infty$ and $E'$ is the sum of all the $h'$-exceptional prime divisors. Let $\psi_\infty':W\to X_\infty'$ be a resolution which resolves the map $X\dashrightarrow Z_m$, and we denote $\psi:W\to X$ the induced morphism. By \cite[Lemma 3.10]{Li22}, for any sufficiently large integer $n$, we have that
$$\kappa\left(F_m,\left(\left(\psi\right)^{-1}_*D+nE_W\right)|_{F_m}\right)=0,$$ 
where $F_m$ is a very general fiber of $W\to Z_m$ and $E_W$ is the sum of all the $\psi$-exceptional prime divisors. Let $E_1$ be the sum of all the $\psi_\infty'$-exceptional prime divisors. Since $E_W=(\psi_\infty')^{-1}_*E'+E_1$ and
$$\psi_*^{-1}D+n\left(\psi_\infty'\right)^{-1}_*E'+(n+m)E_1-\left(\psi_\infty'\right)^*\left(\left(h'\right)_*^{-1}D+nE'\right)\ge0$$
for any sufficiently large integers $n$ and $m$, we see that
$$\kappa\left(F_m,\left(\left(\psi_\infty'\right)^*\left(h_*^{-1}D+nE_1\right)\right)|_{F_m}\right)=0.$$
Note that $Z_\infty$ is birational to $Z_m$, thus there exists a birational morphism from $F_m$ to $F'$. Therefore $\kappa(F',((h')_*^{-1}D+nE')|_{F'})=0$. This finishes the proof.
\end{proof}

\begin{lem}\label{lem:replaceif1}
Notation as in Definition \ref{defn: r cartier r dvisor iitaka fibration}. Suppose that $\phi:X'\to X$ is a birational morphism from a normal projective variety and $D'$ is an $\Rr$-Cartier $\Rr$-divisor on $X'$ such that $D'-\phi^*D\ge0$ and is $\phi$-exceptional. Then 
\begin{enumerate}
    \item possibly replacing $X_\infty$ with a high model, $f_\infty$ is an Iitaka fibration of $D'$, and
    \item any Iitaka fibration of $D'$ is an Iitaka fibration of $D$. 
\end{enumerate}
\end{lem}

\begin{proof}
By our assumption, $H^0(X,\lf mD\rf)=H^0(X_\infty,\lf mD_\infty\rf)$ and $X\dashrightarrow Z_m$ is the map associated with the complete linear system $|\lf mD\rf|$ if and only if $X'\dashrightarrow Z_m$ is the map associated with the complete linear system $|\lf mD'\rf|$.

(1) By Lemma \ref{lem:replaceif0}(2), possibly replacing $X_\infty$ with a high model, we may assume that the map $X_\infty\dashrightarrow X'$ is a morphism and we denote it by $h_0$. It suffices to prove that for any sufficiently large integer $n$,
$$\kappa\left(F,\left((h_0)^{-1}_*D'+nE_0\right)|_F\right)=0,$$ where $E_0$ is the sum of all the $h_0$-exceptional prime divisors. This follows from the fact that 
$$(h_0)^{-1}_*D'+nE_0\le h^{-1}_*D+nE$$ 
for any sufficiently large integer $n$.

(2) Suppose that $f_\infty':X_\infty'\to Z_\infty'$ is an Iitaka fibration of $D'$. Let $W\to X_\infty'$ and $W\to X_\infty$ be a common resolution of $X_\infty$ and $X_\infty'$. By Lemma \ref{lem:replaceif0}(2), $W\to Z_\infty$ is also an Iitaka fibration of $D$. 
It follows that $X_\infty'\to Z_\infty$ is also an Iitaka fibration of $D$. Since $Z_\infty$ is birational to $Z_\infty'$, $X_\infty'\to Z_\infty'$ is an Iitaka fibration of $D$.
\end{proof}

\begin{lem}\label{lem:replaceif2}
Notation as in Definition \ref{defn: r cartier r dvisor iitaka fibration}. Assume that $\psi:X\dashrightarrow X'$ is a $D$-non-negative birational contraction. Then possibly replacing $X_\infty$ with a high model, $f_\infty$ is an Iitaka fibration of $D'$.
\end{lem}
\begin{proof}
By Lemma \ref{lem:replaceif0}(2), possibly replacing $X_\infty$ with a high model, we may assume that the map $h':X_\infty\dashrightarrow X'$ is a morphism. By our assumption, $0\leq h^*D-h'^*D'$ is $h'$-exceptional, and $X'\dashrightarrow Z_m$ is the map associated with the complete linear system $|\lf mD'\rf|$. By Lemma \ref{lem:replaceif0}(1), $f_\infty$ is an Iitaka fibration of $h^*D$. Then by Lemma \ref{lem:replaceif1}(2), $f_\infty$ is also an Iitaka fibration of $D'$.
\end{proof}

\begin{rem}
In the rest of this paper, we will use Lemmas \ref{lem:replaceif0}, \ref{lem:replaceif1} and \ref{lem:replaceif2} frequently to replace $X$ with another birational model without citing.
\end{rem}

\begin{lem}\label{lem: general fiber very general are ok}
Assume that $(X,C)$ is a projective klt pair such that $K_X+C\sim_\Rr0$. Assume that Conjecture \ref{conj: non-vanishing} holds in dimension $\le\dim X$. Suppose that $D$ is a $\Qq$-Cartier $\Qq$-divisor on $X$ such that $\kappa(X,D)\ge0$, and $f_\infty:X_\infty\to Z_\infty$ is an Iitaka fibration of $D$. Then $\kappa(F,h^*D|_F)=0$, where $F$ is a general fiber of $f_\infty$ and $h:X_\infty\to X$ is the induced birational morphism.
\end{lem}
\begin{proof}
By Shokurov type polytopes, we may assume that $C\in\Qq$ and thus $K_X+C\sim_\Qq0$.
Take a positive rational number $\epsilon$ such that $(X_\infty,h^{-1}_*C+(1-\epsilon)E)$ is klt and 
$$E':=K_{X_\infty}+h^{-1}_*C+(1-\epsilon)E-h^*(K_X+C)\ge0,$$
where $E$ is the sum of all the $h$-exceptional prime divisors. Since $\kappa(X,D)\ge0$, one can find a $\Qq$-Cartier $\Qq$-divisor $D'\ge0$ such that $D'\sim_\Qq D$. Let $\epsilon'$ be a positive rational number such that $(X_{\infty},h^{-1}_*C+(1-\epsilon)E+\epsilon'h^*D')$ is still klt. As $f_\infty$ is an Iitaka fibration of $D$ and $E'$ is $h$-exceptional, one can see that
$$\kappa\left(F_v,\left(K_{X_\infty}+h^{-1}_*C+(1-\epsilon)E+\epsilon'h^*D'\right)|_{F_v}\right)=\kappa\left(F_v,\left(E'+\epsilon'h^*D'\right)|_{F_v}\right)=0,$$
where $F_v$ is a very general fiber of $f_\infty$. According to \cite[Corollary 1.4]{HMX18}, we have
$$\kappa\left(F,\left(K_{X_\infty}+h^{-1}_*C+(1-\epsilon)E+\epsilon'h^*D'\right)|_{F}\right)=0.$$
Therefore $\kappa(F,h^*D|_F)=0$.
\end{proof}


\section{Canonical bundle formulas}\label{sec3}
We refer the reader to \cite{BZ16,Bir19,HL21a,HL22a} for the definitions and basic properties for generalized pair (g-pair for short), and we denote by $(X/Z,B+\bM)$ a g-pair throughout this paper. We refer the reader to \cite{Bir19,HL21b,JLX22} for the definition and basic properties of the canonical bundle formula. To sum up, given an lc pair $(X/Z,B)$ and a contraction $\phi:X\to T$ between normal quasi-projective normal varieties over $Z$ such that $K_X+B\sim_{\Rr,T}0$. Then we can find an $\Rr$-divisor $B_T\ge0$ and a nef over $Z$ b-$\Rr$-divisor $\bM_\phi$ on $T$, such that $(T/Z,B_T+\bM_\phi)$ is a glc g-pair, and
$$K_X+B\sim_\Rr\phi^*\left(K_T+B_T+\bM_{\phi,T}\right).$$
Here $B$ (resp., $\bM_\phi$) is called the \emph{discriminant part} (resp., a \emph{moduli part}) of the canonical bundle formula for $(X/Z,B)$ over $T$ which is uniquely determined (resp., determined only up to $\Rr$-linear equivalence). We may also call $\bM_{\phi,T}$ the moduli part of the canonical bundle formula for $(X/Z,B)$ over $T$. Moreover, if $(X/Z,B)$ is klt, then $(T/Z,B_T+\bM_\phi)$ is gklt.

Here we emphase that there are many choices of $\bM_\phi$, some of which could behave badly. But we can always choose one with the required properties in the following results.

For convenience, we say that two g-pairs $(X/Z,B+\bM)$ and $(X'/Z,B'+\bM')$ are \emph{crepant} if $X$ is birational to $X'$, $\bM=\bM'$, and $p^*(K_X+B+\bM_X)=q^*(K_{X'}+B'+\bM'_{X'})$ for some common resolution $p:W\to X$ and $q:W\to X'$. We also call $(X'/Z,B'+\bM)$ a \emph{crepant model} of $(X/Z,B+\bM)$.

\begin{prop}\label{prop:cbfindex}
Let $d$ be a positive integer and $\Phi\subset[0,1]$ a DCC set. Assume that Conjectures \ref{conj: boundedness and existence of n complement nft} and \ref{conj: exist gmm} hold in dimension $d$. Then there exist a positive integer $p$ and a DCC set $\Phi'$ depending only on $d$ and $\Phi$ satisfying the following.

Assume that $(X/Z,B)$ is an lc pair of dimension $d$ and $\phi:X\to T$ is a contraction over $Z$ such that $\dim T>0$, $B\in\Phi$, $K_X+B\sim_{\Rr,T}0$, and $K_X+B\sim_{\mathbb Q,T}0$ over the generic point of $T$. Then we can choose a moduli part $\bM_\phi$ of the canonical bundle formula for $(X/Z,B)$ over $T$ such that $B_T\in\Phi'$, $p\Mm_\phi$ is $\bb$-Cartier, and
$$p(K_X+B)\sim p\phi^*(K_T+B_T+\Mm_{\phi,T}),$$
where $B_T$ is the discriminant part of the canonical bundle formula for $(X/Z,B)$ over $T$. 

Moreover, if $\Phi$ is a hyperstandard set, then $\Phi'$ is a hyperstandard set.
\end{prop}
The proof is similar to the proof of \cite[Proposition 6.3]{Bir19}. For the convenience of the reader, we give a proof here. We also remark that the moreover part of the proposition will not be used in this paper, but it is useful in some other situations (cf. \cite{CHX22}).

\begin{proof}

\noindent{\bf Step 1}. In this step, we construct $p$ and make a choice of $\bM_{\phi,T}$. Note that here we only make a choice of $\bM_{\phi,T}$ rather than $\bM_\phi$.

By \cite[Theorems 1.8 and 8.25]{HLS19}, there exist a positive integer $p$ and a finite set $\Ii_0\subset (0,1]$ depending only on $d$ and $\Phi$, such that for any $t\in T$, $(X/T\ni t,B)$ has a $(p,\Ii_0)$-decomposable $\Rr$-complement $(X/T\ni t,B+G)$ for some $\Rr$-Cartier $\Rr$-divisor $G\ge0$, and moreover if $B\in\Ii\cap\Qq$, then $(X/T\ni t,B)$ has a monotonic $p$-complement. In particular, $G\sim_{\Rr}0$ over a neighborhood of $t$, and hence $G$ is vertical over $T$. Since $K_X+B\sim_{\mathbb Q,T}0$ over the generic point $\eta_T$ of $T$, $p(K_X+B)\sim0$ over a neighborhood of $\eta_T$, and there exists $\alpha\in \mathcal{K}(X)$ such that $pL:=p(K_X+B)+(\alpha)$ is zero near $\eta_T$. In particular, $pL\sim p(K_X+B)\sim_{\Rr,T}0$ and $L$ is vertical over $T$. By \cite[Lemma 2.11]{Li20}, we may find an $\Rr$-Cartier $\Rr$-divisor $L_T$ on $T$ such that $L=\phi^*L_T$. Let $B_T$ be the discriminant part of the canonical bundle formula of $(X,B)$ over $T$, and $\bM_{\phi,T}:=L_T-K_T-B_T.$ Then
$$p(K_X+B)\sim pL=p\phi^*L_T=p\phi^*(K_T+B_T+\bM_{\phi,T}).$$

\medskip

\noindent{\bf Step 2}. In this step, we show that we can reduce to the case $\dim T=1$ to show the existence of $\Phi'$ and prove that $p\bM_{\phi,T}$ is integral. 

Assume $\dim T>1$. Let $H$ be a general hyperplane section of $T$, $G:=\phi^*H$, and $g:G\to H$ the induced morphism. We may write $K_{G}+B_{G}=(K_{X}+G+B)|_{G}.$ It is clear that $(G,B_G)$ is an lc pair with $K_G+B_G\sim_{\Qq,H}0$. Suppose that $B_{H}$ is the discriminant part of the canonical bundle formula for $(G, B_{G})$ over $H$. Note that as $G$ is a general member of a free linear system, every lc center $S_G$ of $(G,B_G)$ is a component of $S_0\cap G$ for some lc center $S_0$ of $(X,B)$.

We claim that $\mult_{D} B_{T}=\mult_{C} B_{H}$ for any prime divisor $D$ on $T$ and any component $C$ of $D \cap H$. Indeed, let $s_D$ be the lc threshold of $\phi^{*}D$ with respect to $(X,B)$ over the generic point of $D$. Then there is an lc center $F$ of $(X, B+s_D\phi^{*}D)$ such that $\phi(F)=D$. Note that $F$ is also an lc center of $(X, B+G+s_D\phi^{*}D)$ as $G$ is general. Hence $F\cap G$ is an lc center of $ (G, B_{G}+s_Dg^{*}C)$, by inversion of adjunction \cite{Kawakita07}. Moreover, as $\phi(F\cap G)=C$, we see that $s_D$ is the lc threshold of $g^{*} C$ with respect to $(G,B_{G})$ over the generic point of $C$, and the claim holds. Since $\Phi$ is a DCC set, there is a DCC set $\Phi_1$ depending only on $\Phi$ such that $B_{G}\in \Phi_1$ (cf. \cite[Corollary 16.7]{Kol92}). If $\Phi$ is a hyperstandard set, then we can take $\Phi_1$ to be a hyperstandard set by \cite[Lemma 3.3]{Bir19}. In both case, by induction, there is a DCC set (resp. hyperstandard set) $\Phi_1'$ such that $B_{Z} \in\Phi_1'$.

Pick a general $H' \sim H$ and let $K_{H}:= (K_{T}+H')|_{H} .$ Note that the restriction is well defined as $H$ is a general hyperplane section and $K_{H}$ is determined as a Weil divisor, although $K_{T}$ may not be $\Qq$-Cartier. Let
$$\bM_{g,H}:= (L_{T}+H')|_{H}- (K_{H}+B_{H} ).$$
Then $B_{H}+\bM_{g,H}=(B_{T}+\bM_{\phi,T})|_{H}$, and
$$p (K_{G}+B_{G} ) \sim p(L+G)|_{G} \sim pg^{*}(L_{T}+H')|_{H} \sim p g^{*}(K_{H}+B_{H}+\bM_{g,H}).$$ 
Hence $\bM_{g,H}$ is the moduli part of the canonical bundle formula for $(G,B_{G})$ over $H$, and
$$\mult_{C} (B_{H}+\bM_{g,H} )=\mult_{D} (B_{T}+\bM_{\phi,T} )$$ 
which implies that $\mult_{C} \bM_{g,H}=\mult_{D} \bM_{\phi,T}$ as $\mult_{C} B_{H}=\mult_{D} B_{T}.$ Therefore $p\mult_{D}\bM_{\phi,T}$ is integral if and only if $p\mult_{C}\bM_{g,H}$ is integral. Repeating the process we may finish this step.
\medskip

\noindent{\bf Step 3}. In this step we show the existence of $\Phi'$ and that $p\bM_{\phi,T}$ is integral. Note that by {\bf Step 2}, we may assume that $\dim T=1$.

\smallskip

\noindent{\bf Step 3.1}. We construct the set $\Phi'$. If $\Phi$ is a hyperstandard set which is not a hyperstandard set, then by \cite[Theorem 1.1]{HMX14}, $B_T\in\Phi'$ for some DCC set $\Phi'$ which only depends on $d$ and $\Phi$. If $\Phi=\Phi(\fR)$ is a hyperstandard set, then we show $B_T\in\Phi':=\Phi(\fR')$, where
$$\fR':=\left\{r,\frac{r}{l_1}-\frac{l_2}{p}\mid r\in\fR,\, l_1\,, l_2\in\Zz_{>0}\right\}\cap[0,1]$$
is a finite set of rational numbers.

\begin{claim}\label{claim:lctcmpt}
Suppose that $t\in T$ is a closed point. Then $(X/T\ni t,B+s\phi^*t)$ is a monotonic $p$-complement of $(X/T\ni t,B)$, where $s:=\lct(X/T\ni t,B;\phi^*t)$. In particular, $B+s\phi^*t$ is a $\Qq$-divisor over a neighborhood of $t$.
\end{claim}
\begin{proof}[Proof of Claim \ref{claim:lctcmpt}]
By Lemma \ref{lem: dim 1 base q divisor}, $B_t:=B+s\phi^*t$ is a $\Qq$-divisor over a neighborhood of $t$. Possibly shrinking $T$ around $t$, we may assume that $B_t\in\Qq$. Let $(X',B_t')$ be a $\Qq$-factorial dlt modification of $(X,B_t)$. Then $\lf B_t'\rf$ has a component mapping to $t$ and $K_{X'}+B_t'\sim_{\Qq,T}0$. There exists a $\Qq$-divisor $B'$ on $X'$ such that $B'\in\Phi\cap\Qq$, $\lf B'\rf=\lf B_t'\rf$, and $B_{X'}\le B'\le B_t'$, where $B_{X'}$ is the strict transform of $B$ on $X'$. By assumption, $(X'/T\ni t,B')$ has a monotonic $p$-complement $(X'/T\ni t,B'^{+})$. Let $B^+$ be the strict transform of $B'^{+}$ on $X$. Then $(X/T\ni t,B^{+})$ is a monotonic $p$-complement of $(X/T\ni t,B)$. Since $B^+-B\ge0$ and $B^+-B\sim_{\Qq}0$ over a neighborhood of $t$, $B^+-B$ is vertical over $T$. Moreover, as $B'^{+}\ge B'$, $\lf B'^{+}\rf$ has a component mapping to $t$, and thus $(X,B^+)$ has an lc center mapping to $t$. Therefore $B^+-B=s\phi^*t$ over $t$. The claim holds.
\end{proof}

Pick a closed point $t\in T$. By Claim \ref{claim:lctcmpt}, $(X/T\ni t,B^+:=B+s\phi^*t)$ is a $p$-complement of $(X/T\ni t,B)$, where $s:=\lct(X/T\ni t,B;\phi^*t)$. For any component $S$ of $\phi^*t$, let $b:=\mult_SB,\, b^+:=\mult_SB^+$ and $m:=\mult_S\phi^*t$. Then $b^+=b+sm$ and thus $s=\frac{b^+-b}{m}$. Since $b\in\Phi$, we may write $b=1-r/l$ for some $r\in\fR$ and $l\in\Zz_{>0}$. In particular, 
$$s=\frac{b^+-1+r/l}{m}=\frac{1}{m}\left(\frac{r}{l}-\left(1-b^+\right)\right).$$
Then $\mult_tB_T=1-s\in\Phi',$ as $b^+\in\frac{1}{p}\Zz\cap[0,1]$.

\smallskip

\noindent{\bf Step 3.2}. We show that $p\bM_{\phi,T}$ is integral. 

We may assume that $p(K_X+B)\sim0$ over some non-empty open subset $U_0\subseteq T$ such that $\Supp B_T\subseteq T\setminus U_0$. Let 
$$\Theta:=B+\sum_{t\in T\setminus U_0}s_{t}\phi^*t,$$
where $s_t:=\lct(X/T\ni t,B;\phi^*t)$. Let $\Theta_T$ be the discriminant part of the canonical bundle formula for $(X,\Theta)$ over $T$. Then
$$\Theta_T=B_T+\sum_{t\in T\setminus U_0}s_{t}t$$ 
which is a reduced divisor. Moreover, $(X/T\ni t,\Theta)$ is a $p$-complement of $(X/T\ni t,B)$ for every $t\in T\setminus U_0$ by Claim \ref{claim:lctcmpt}. Hence $p(K_X+\Theta)\sim_T0$ by Lemma \ref{lem:trivialoverbase}. Since 
\begin{align*}
\ \ \ \ p(K_X+\Theta)&=p(K_X+B)+p(\Theta-B)\\& \sim p\phi^* (K_T+B_T+\bM_{\phi,T} )+p \phi^*(\Theta_T-B_T)\\&=p\phi^*(K_T+\Theta_T+\bM_{\phi,T}),
\end{align*}
$p(K_T+\Theta_T+\bM_{\phi,T})$ is Cartier. It follows that $p\bM_{\phi,T}$ is integral as $K_T+\Theta_T$ is reduced.

\medskip

\noindent{\bf Step 4}. In this step, we finish the proof by showing that $p\bM_\phi$ is b-Cartier and nef. Note that in this step we do not assume $\dim T=1$.

According to \cite[Theorem 3.6]{Bir19}, we only need to show that $p\bM_\phi$ is b-Cartier. Let $T'\to T$ be a high resolution and $Y \to X$ a log resolution of $(X, B)$ such that $Y\to T'$ is a morphism. Let $U \subseteq T$ be a non-empty open subset such that $U'\to U$ is an isomorphism where $U'\subseteq T'$ is the inverse image of $U$. Let $B_Y$ be the sum of the strict transform of $B$ and the reduced exceptional divisor of $Y  \to X$ but with all the components mapping outside $U$ removed. In particular, the generic point of any lc center of $(Y, B_Y)$ maps into $U$. We may run an MMP on $K_{Y}+B_Y$ over $X\times_{T} T'$ with scaling of some ample divisor. By \cite[Theorem 1.9]{Bir12}, the MMP terminates over $U'$. In fact, we reach a model $X'$ such that over $U'$ the pair $(X', B')$ is a $\Qq$-factorial dlt modification of $(X, B)$, where $B'$ is the strict transform of $B_Y$ on $X'$. Hence $K_{X'}+B'\sim_{\Qq} 0$ over $U'$. Now by \cite[Theorem 1.1]{HX13} (see also \cite[Theorem 1.1]{Bir12}), we can run an MMP on $K_{X'}+B'$ over $T'$ which terminates with a good minimal model $X''$ over $T'$ as the generic point of every lc center of $(X', B')$ is mapped into $U'$. Let $B''$ be the strict transform of $B'$ on $X''$. Then $K_{X''}+B''$ is semi-ample over $T'$. 

\begin{center}
  \begin{tikzcd}[column sep = 2em, row sep = 2.5em]
     X   \arrow[dr, "\phi",swap] && \arrow[ll, ""] Y\arrow[dr, ""] \arrow[rr, "",dashed]&& X' \arrow[dl, ""]\arrow[rr, "",dashed]&& X'' \arrow[dl, "\phi''"]\\&
   T &&  \arrow[ll, ""]  T'   &&\arrow[ll, ""] T'' \\&
   U \arrow[hook]{u}&&  \arrow[ll, "\cong" swap] U'\arrow[hook]{u}   &&\arrow[ll, "\cong" swap]U''\arrow[hook]{u}
  \end{tikzcd}
\end{center}

Let $\phi'': X''\to T'' $ be the contraction defined by $K_{X''}+B''$ over $T'$. Note that $T''\to T'$ is birational as $K_{X'}+B'\sim_{\Qq} 0$ over $U'$. Assume that $(T'',B_{T''}''+\bM_\phi)$ is the crepant model of $(T,B_T+\bM_\phi)$. Then $L_{T''}=K_{T''}+B_{T''}''+\bM_{\phi,T''}$, where $L_{T''}$ is the pullback of $L_T$ on $T''$. Let $f:W\to X$ and $f'':W\to X''$ be a common resolution, $K_{X''}+\Delta'':=f''_*f^*(K_X+B)$. Since $K_{X}+B$ and $K_{X''}+B''$ are crepant over $U$, we see that $B''-\Delta''$ is vertical over $T''$. Note that $B''-\Delta''\sim_{\Rr,T''}0$, $B''-\Delta''=(\phi'')^*P_{T''}$ for some $\Rr$-Cartier $\Rr$-divisor $P_{T''}$ on $T''$ by Lemma \ref{lem:efftri}. Denote by $B_{T''}$ the discriminant part of the canonical bundle formula for $(X'', B'')$ over $T''$. Then $B_{T''}=B_{T''}''+P_{T''}$ and 
$$
\begin{aligned}
p\left(K_{X''}+B''\right) &=p\left(K_{X''}+\Delta''+B''-\Delta''\right) \sim p \left(f''_*f^*L+B''-\Delta''\right)\\&=p\left(\phi''\right)^{*} \left(L_{T''}+P_{T''} \right)=p \left(\phi''\right)^{*} \left(K_{T''}+B_{T''}+\bM_{\phi,T''} \right).
\end{aligned}
$$
Now by {\bf Step 3.2}, $p\bM_{\phi,T''}$ is an integral divisor, hence $p\bM_{\phi,T'}$ is integral. As $T'$ is smooth, we conclude that $p\bM_\phi$ is b-Cartier.
\end{proof}

\begin{lem}\label{lem:normlzt}
Assume that $X\to Z$ is a proper morphism of varieties. Let $X^{\nu}$ be the normalization of $X$. Then $X^{\nu}_z$ is the normalization of $X_z$ for any general point $z\in Z$, where $X_z$ (resp., $X^{\nu}_z$) is the fiber over $z$. In particular, if $X_z$ is normal, then $X^{\nu}_z$ is isomorphic to $X_z$.
\end{lem}
\begin{proof}
By assumption, $X^\nu\to X$ is birational and finite. Thus $X^\nu_z\to X_z$ is birational and finite. Let $\tilde X_z$ be the normalization of $X_z$. The universal property of the normalization implies that there is a morphism $\tilde X_z\to X^\nu_z$. Moreover, the morphism $\tilde X_z\to X^\nu_z$ is birational and finite. Since both $\tilde X_z$ and $X^\nu_z$ are normal, we see that $\tilde X_z$ is isomorphic to $X^\nu_z$. 
\end{proof}

\begin{prop}\label{prop:cbfindex2}
Let $d,b,\beta$ be three positive integers and $\Ii\subset [0,1]\cap\Qq$ a finite set. Then there exist a positive integer $p$ and a DCC set $\Ii'\subseteq[0,1]$ depending only on $d,b,\beta$ and $\Ii$ such that $\bar{\Ii'}\subseteq[0,1]\cap\Qq$ and satisfying the following. 

Assume that $(X,B)$ is an lc pair of dimension $d$, $\phi:X\to T$ is a contraction between quasi-projective normal varieties, and $F$ is a general fiber of $\phi$. Suppose that
\begin{enumerate}
    \item $B\in\Ii$, and $K_X+B\sim_{\Qq,T}0$,
    \item $\dim T>0$ and $(X,B)$ is klt over the generic point of $T$,
    \item $b$ is the non-vanishing order of $(K_X+B)|_F$, and
    \item $\beta:=\dim H^{\dim F}(\tilde F,\mathbb C)$, where $\tilde F$ is a smooth model of the cover of $F$ associated to the unique divisor of $|b(K_X+B)|_F|$.
\end{enumerate}
Then we can choose a moduli part $\bM_\phi$ of the canonical bundle formula for $(X,B)$ over $T$ such that $B_{T}\in\Ii'$, $p\bM_\phi$ is b-Cartier, and
$$p\left(K_X+B\right)\sim p\phi^*\left(K_{T}+B_{T}+\bM_{\phi,T}\right),$$
where $B_{T}$ is the discriminant part of the canonical bundle formula for $(X,B)$ over $T$.
\end{prop}
\begin{proof}
Let $B_{T}$ be the discriminant part of the canonical bundle formula for $(X,B)$ over $T$. Then, by \cite[Theorem 1.11]{HMX14}, $B_T\in\Ii'$ for some DCC set $\Ii'$ which only depends on $d$ and $\Ii$ such that $\bar{\Ii'}\subset[0,1]\cap\Qq$. It suffices to prove the existence of $p$ with the required properties.

Let $T'\to T$ be a resolution, $X'$ the normalization of the main component of $X\times_TT'$, and denote by $\phi':X'\to T'$. Let $(X',B')$ be the crepant model of $(X,B)$. Then $K_{X'}+B'\sim_{\Qq,T'}0$ and $(X',B')$ is sub-lc and is klt over the generic point of $T'$. By Lemma \ref{lem:normlzt}, the general fiber $F'$ of $\phi'$ is isomorphic to $F$. Therefore it holds that
\begin{itemize}
  \item $b$ is the non-vanishing order of $(K_{X'}+B')|_{F'}$, and
  \item $\tilde F$ is a smooth model of the cover of $F'$ associated to the unique divisor of $|b(K_{X'}+B')|_{F'}|$.
\end{itemize}
According to \cite[Theorem 5.1]{Flo14}, one can find a positive integer $p_0$ depending only on $\beta$ and a choice of a moduli part $\bM_\phi$ of the canonical bundle formula for $(X',B')$ over $T'$ such that $p_0\bM_{\phi,T'}$ is integral. In particular, $p_0\bM_{\phi,T'}$ is Cartier and thus $p_0\bM_\phi$ is b-Cartier. Since $b$ is the non-vanishing order of $(K_X+B)|_F$, by the choice of $\bM_\phi$ again, we see that
$$b\left(K_X+B\right)\sim b\phi^*\left(K_{T}+B_{T}+\bM_{\phi,T}\right).$$ 
Therefore we may conclude that $p:=bp_0$ and $\Ii'$ have the required properties.
\end{proof}

\section{Non-vanishing orders, middle Betti numbers, and complements}\label{sec6}

In this section, we prove Theorem \ref{thm: betti to complement intro}. More generally, we proof Theorem \ref{thm: betti to complement intro1}. 

\begin{thm}\label{thm: betti to complement intro1}
Let $d,b,\beta$ be three positive integers and $\Ii\subset [0,1]\cap\Qq$ a DCC set. Assume that the good minimal model conjecture (Conjecture \ref{conj: exist gmm}) holds in dimension $d$. Then there exists a positive integer $n$ depending only on $d,b,\beta$ and $\Ii$ satisfying the following. 

Assume that $(X,B)$ is a $\Qq$-factorial projective pair of dimension $d$, $f_\infty:X_\infty\to Z_\infty$ is an Iitaka fibration of $-(K_X+B)$, $h:X_\infty\to X$ is the induced birational morphism, and $F$ is a general fiber of $f_\infty$. Suppose that
\begin{enumerate}
    \item $B\in\Ii$,
    \item $-(K_X+B)$ is nef,
    \item $(X,B)$ has a klt $\Rr$-complement,
    \item $\kappa(X,-(K_X+B))\ge0$, and $b$ is the non-vanishing order of $-h^*(K_X+B)|_F$, i.e.,
    $$b=\min\left\{a\in\Zz_{>0}\mid |-ah^*(K_X+B)|_F|\not=\emptyset\right\},$$
and
    \item $\beta:=\dim H^{\dim F}(\tilde F,\mathbb C)$, where $\tilde F$ is a smooth model of the cover of $F$ associated to the unique divisor of $|-bh^*(K_X+B)|_F|$. 
\end{enumerate}
Then $(X,B)$ has an $n$-complement.
\end{thm}

We remark here that $\kappa(F,-h^*(K_X+B)|_F)=0$ by Lemma \ref{lem: general fiber very general are ok}. We first show the existence of $(n_0,\Ii_0)$-decomposable $\Rr$-complements under the assumption of Theorem \ref{thm: betti to complement intro1}.

\begin{thm}\label{thm: betti to decomposable complement}
Notation and assumptions as in Theorem \ref{thm: betti to complement intro1}. Then $(X,B)$ has an $(n_0,\Ii_0)$-decomposable $\Rr$-complement, where $n_0$ is a positive integer and $\Ii_0\subset[0,1]$ is a finite set depending only on $d,b,\beta$ and $\Ii$.
\end{thm}

\begin{proof}
\noindent{\bf Step 1}. In this step, we show that $-(K_X+B)$ is semi-ample.

Let $(X,C)$ be a klt $\Rr$-complement of $(X,B)$ for some effective $\Rr$-divisor $C\ge B$. Pick some positive real number $\epsilon_0$ such that $(X,C+\epsilon_0(C-B))$ is klt. Since $K_X+C+\epsilon_0(C-B)\sim_\Rr -\epsilon_0(K_X+B)$ is nef, we see that $K_X+C+\epsilon_0(C-B)$ is semi-ample by \cite[Lemma 2.9.1]{HMX18} and hence $-(K_X+B)$ is also semi-ample. Denote by $f:X\to Z$ the ample model of $-(K_X+B)$ and $F_0$ a general fiber of $f$. Moreover, by \cite[Proposition 3.21]{Li22}, $f_\infty$ is birational to $f$, so there exists a naturally induced morphism $h: F\rightarrow F_0$ such that $b(K_X+B)|_{F_0}\sim0$.

If $\dim Z=0$, then $b(K_X+B)\sim0$ and thus $(X,B)$ has a $(b,\{1\})$-decomposable $\Rr$-complement. Therefore, we may assume that $\dim Z>0$.

\medskip

\noindent{\bf Step 2}. In this step, we construct finite sets $\Ii_0\subset(0,1]$ and $\Ii_2'\subset[0,1]\cap\Qq$ depending only on $d,b$ and $\Ii$.

Let $\alpha$ be a positive real number such that 
$$\alpha<\min\left\{|\gamma-\frac{i}{b}|>0\mid \gamma\in\bar\Ii,i\in\Zz_{\ge0}\right\}.$$ 
By \cite[Theorem 5.20]{HLS19}, there exist a finite set $\Ii_1$ depending only on $d,\alpha$ and $\Ii$, and an $\Rr$-divisor $\bar B$ on $X$, such that
\begin{itemize}
    \item $\bar B\in\Ii_1$,
    \item $\alpha\Supp B\ge\bar B-B\ge0$, and
    \item $(X,\bar B)$ is $\Rr$-complementary.
\end{itemize}
In particular, for any component $S$ of $B$ such that $b\mult_SB\in\Zz$, $\mult_S\bar B=\mult_SB$. Hence $(\bar B)^h=B^h$. By \cite[Theorem 5.16]{HLS19}, there exist a point $\bm{v}_0:=(v^0_1,\dots,v^0_m)$ and an open subset $V$ of the rational envelope of $\bm{v}_0$ (that is the smallest affine subspace containing $\bm{v}_0$ which is defined over $\Qq$) in $\Rr^m$ depending only on $d$ and $\Ii_1$, and Weil divisors $B_1,\dots,B_m\geq 0$ on $X$, such that $B(\bm{v}_0)=\bar B$, and $(X,B(\bm{v}))$ is $\Rr$-complementary for any $\bm{v}\in V$, where $B(\bm{v}):=\sum_{i=1}^mv_iB_i$ for any $\bm{v}:=(v_1,\dots,v_m)\in\Rr^m$. Moreover, possibly replacing $V$, we may assume that 
$$B(\bm{v})\ge D:=\bar B-B\text{, and }|B(\bm{v})-\bar B|<\alpha$$ 
for any $\bm{v}\in V$. We remark that
$(B(\bm{v}))^h=(\bar B)^h=B^h$ for any $\bm{v}\in V$.

Pick points $\bm{v}_1,\dots,\bm{v}_k\in V\cap\Qq^m$, such that $\bm{v}_0$ is in the interior of the convex hull spanned by $\bm{v}_1,\dots,\bm{v}_k$. For any integer $1\le i\le k$, set 
$$B^{(i)}:=B(\bm{v}_i)-D\text{ and }\bar B^{(i)}:=B(\bm{v}_i).$$ 
There exist finite sets $\Ii_0:=\{a_1,\dots,a_k\}\cup\{1\}\subset (0,1]$ and $\Ii_2'\subset[0,1]\cap\Qq$ such that $\bar B^{(i)}\in\Ii_2'$ for any integer $1\le i\le k$, $\sum_{i=1}^ka_i=1$, and $\sum_{i=1}^ka_i\bm{v}_i=\bm{v}_0$. In particular, $\sum_{i=1}^k a_iB^{(i)}=B$.

\medskip

\noindent{\bf Step 3}. In this step, we run an MMP.

Let $(X,\bar B^{(i)}+\bar G_i)$ be an $\Rr$-complement of $(X,\bar B^{(i)})$ for some effective $\Rr$-Cartier $\Rr$-divisor $\bar G_i$ for any integer $1\le i\le k$. Pick a positive real number $\epsilon$ such that $(X,C+\epsilon\bar G_i)$ is klt. Then we may run an MMP on $K_X+C+\epsilon\bar G_i$ over $Z$, which terminates with a model $W_i$ over $Z$ such that $K_{W_i}+C_{W_i}+\epsilon\bar G_{W_i}$ is semi-ample over $Z$, where $C_{W_i}$ and $\bar G_{W_i}$ are the strict transforms of $C$ and $\bar G_i$ on $W_i$ respectively. This MMP is also an MMP on $-(K_X+\bar{B}^{(i)})$ over $Z$, hence $-(K_{W_i}+\bar{B}_{W_i}^{(i)})$ is semi-ample over $Z$. Let $g_i: W_i\rightarrow T_i$ be the ample model of $-(K_{W_i}+\bar{B}_{W_i}^{(i)})$ over $Z$. Recall that by construction and by Lemma \ref{lem:horareQ}, $B-\bar{B}^{(i)}$ is vertical over $Z$. As $-(K_X+\bar{B}^{(i)})\sim_{\Rr,Z} B-\bar{B}^{(i)}$, one can see that $F_0$ is isomorphic to a general fiber of $W_i\to Z$ by Lemma \ref{lem:gfiso}, and $T_i$ is birational to $Z$.

\medskip

\noindent{\bf Step 4}. In this step, we construct a positive integer $p$, and a hyperstandard set $\Ii_2$ depending only on $d,b,\beta$ and $\Ii_2'$, and make a choice of a moduli part $\bM_i$ of the canonical bundle formula for $(W_i,\bar B_{W_i}^{(i)})$ over $T_i$ such that
$$p\left(K_{W_i}+\bar B_{W_i}^{(i)}\right)\sim pg_i^*\left(K_{T_i}+B_{T_i}+\Mm_{i,T_i}\right),$$
$B_{T_i}\in\Ii_2$, and $p\Mm_i$ is b-Cartier, where $B_{T_i}$ is the discriminant part of the canonical bundle formula for $(W_i,\bar B_{W_i}^{(i)})$ over $T_i$.

To see this, we first claim that 
\begin{itemize}
  \item[(i)] $b$ is the non-vanishing order of $-(K_{W_i}+\bar B_{W_i}^{(i)})|_{F_i'}$, and $b(K_{W_i}+\bar B_{W_i}^{(i)})|_{F_i'}\sim0$, and
  \item[(ii)] $\tilde F$ is a smooth model of the cover of $F_i'$ associated to the unique divisor of $|b(K_{W_i}+\bar B_{W_i}^{(i)})|_{F_i'}|$, where $F_i'$ is a general fiber of $g_i$.
\end{itemize}
Assume that $b_0$ is the non-vanishing order of $(K_X+B)|_{F_0}\sim_\Qq0$. Then $b_0(K_X+B)|_{F_0}\sim0$ and $b_0\le b$. Since $h_0^*((K_X+B)|_{F_0})=h^*(K_X+B)|_F$, we have $b_0=b$ and $bh^*(K_X+B)|_F\sim0$. Therefore $b$ is also the non-vanishing order of $(K_X+B)|_{F_0}$. Note that 
$$\Spec\left(\bigoplus_{i=0}^{b-1}\Oo_F\left(\lf ih^*(K_X+B)|_{F}\rf\right)\right)\to F\ \ (\text{resp., } \Spec\left(\bigoplus_{i=0}^{b-1}\Oo_{F_0}\left(\lf i(K_X+B)|_{F_0}\rf\right)\right)\to F_0)$$
is the cover associated to the unique element of $\lf h^*(K_X+B)|_F\rf$, see \cite[\S2.3]{Kol13}. By \cite[II Lemma 2.11]{Nak04}, there is a natural isomorphism 
$$(h_0)_*\Oo_F\left(\lf i(h^*(K_X+B)|_{F})\rf\right)\to \Oo_{F_0}\left(\lf i((K_X+B)|_{F_0})\rf\right).$$
Thus $\tilde F$ is a smooth model of the cover of $F_0$ associated to the unique divisor of $|b(K_{X}+B)|_{F_0}|$. Furthermore, by our construction, $X$ is isomorphic to $W_i$ and $B=\bar{B}^{(i)}$ over the generic point of $Z$, and $F_i'$ is isomorphic to $F_0$ (see Lemma \ref{lem:gfiso}). Thus the claim holds.

Since $(X,B)$ has a klt $\Rr$-complement, $(W_i,B_{W_i})$ has a klt $\Rr$-complement, where $B_{W_i}$ is the strict transform of $B$ on $W_i$. Moreover, as $B_{W_i}=\bar B_{W_i}^{(i)}$ over the generic point of $T_i$, we can see that $(W_i,\bar B_{W_i}^{(i)})$ is klt over the generic point of $T_i$. Recall that $\Ii_2'\subset[0,1]\cap\Qq$ is a finite set and $\bar B_{W_i}^{(i)}\in\Ii_2'$. By Proposition \ref{prop:cbfindex2}, we may find a positive integer $p$, and a DCC set $\Ii_2$ depending only on $d,b,\beta$ and $\Ii_2'$ such that $\bar \Ii_2\subset[0,1]\cap\Qq$, and make a choice of a moduli part $\bM_i$ of the canonical bundle formula for $(W_i,\bar B_{W_i}^{(i)})$ over $T_i$, such that
$$p\left(K_{W_i}+\bar B_{W_i}^{(i)}\right)\sim pg_i^*\left(K_{T_i}+B_{T_i}+\Mm_{i,T_i}\right),$$
$B_{T_i}\in\Ii_2$, and $p\Mm_i$ is b-Cartier, where $B_{T_i}$ is the discriminant part of the canonical bundle formula for $(W_i,\bar B_{W_i}^{(i)})$ over $T_i$. Note that $(T_i,B_{T_i}+\bM_i)$ is glc.

\medskip

\noindent{\bf Step 5}. In this step, we find an integer $n_0$ which only depends on $d,bp$ and $\Ii_2$ satisfying our requirements and thus finish the proof.

We first show that $T_i$ is of Fano type. Indeed, according to \cite[Theorem 0.2]{Amb05}, there exist klt pairs $(Z,\Delta)$ and $(T_i,\Delta_{T_i})$ such that 
$$K_X+B\sim_\Rr f^*(K_Z+\Delta)\text{ and }K_{W_i}+B_{W_i}\sim_\Rr g_i^*(K_{T_i}+\Delta_{T_i}).$$
In particular, $K_{T_i}+\Delta_{T_i}\sim_\Rr \tau_i^*(K_Z+\Delta)$ where $\tau_i$ denotes the induced morphism $T_i\to Z$. Since $Z$ is the ample model of $-(K_X+B)$, $-(K_Z+\Delta)$ is ample. It follows that $-(K_{T_i}+\Delta_{T_i})$ is big and nef as $\tau_i$ is a birational morphism. Thus $T_i$ is of Fano type.

Since $(W_i,\bar B_{W_i}^{(i)})$ is $\Rr$-complementary, $(T_i,B_{T_i}+\bM_i)$ is $\Rr$-complementary, that is there is an $\Rr$-divisor $P_i\ge0$ such that $(T_i,B_{T_i}+P_i+\Mm_i)$ is glc and $K_{T_i}+B_{T_i}+P_i+\Mm_{i,T_i}\sim_\Rr 0$. As $T_i$ is of Fano type, by \cite[Theorem 1.10]{Bir19} (see also \cite[Theorem 1.3]{Che20}), there exists a positive integer $n_0$ divisible by $bp$ depending only on $d,bp$ and $\Ii_2$, and a $\Qq$-divisor $B_{T_i}^+\geq B_{T_i}$ on $T_i$, such that $(T_i,B_{T_i}^++\Mm_i)$ is glc and $n_0(K_{T_i}+B_{T_i}^++\Mm_{i,T_i})\sim 0$.

It is enough to show that $(W_i,\bar{B}_{W_i}^{(i),+}:=\bar B_{W_i}^{(i)}+g_i^*(B_{T_i}^+-B_{T_i}))$ is a monotonic $n_0$-complement of $(W_i,\bar{B}_{W_i}^{(i)})$. Indeed, by \cite[Corollary 7.18(1)]{PS09}, $(W_i,\bar{B}_{W_i}^{(i),+})$ is lc. Since
\begin{align*}
    &\ \ \ \ n_0\left(K_{W_i}+\bar{B}_{W_i}^{(i),+}\right)=n_0\left(K_{W_i}+\bar B_{W_i}^{(i)}\right)+n_0g_i^*\left(B_{T_i}^+-B_{T_i}\right)\\
    &\sim n_0g_i^*\left(K_{T_i}+B_{T_i}+\Mm_{i,T_i}\right)+n_0g_i^*\left(B_{T_i}^+-B_{T_i}\right)
    \\&=n_0g_i^*\left(K_{T_i}+B_{T_i}^++\Mm_{i,T_i}\right)\sim 0,
\end{align*}
one can see that $(W_i,\bar{B}_{W_i}^{(i),+})$ is a monotonic $n_0$-complement of $(W_i,\bar{B}_{W_i}^{(i)})$. Remark that by {\bf Step 3}, $X\dashrightarrow W_i$ is $-(K_X+\bar{B}^{(i)})$-negative. Therefore $(X,\bar{B}^{(i)})$ also has a monotonic $n_0$-complement $(X,\bar{B}^{(i),+})$. This immediately implies that $(X,\sum_{i=1}^ka_i\bar{B}^{(i),+})$ is an $(n_0,\Ii_0)$-decomposable $\Rr$-complement of $(X,B)$. We may finish the proof.
\end{proof}

\begin{proof}[Proof of Theorem \ref{thm: betti to complement intro1}]
By Theorem \ref{thm: betti to decomposable complement}, there exist a positive integer $n_0$ and a finite set $\Ii_0\subset (0,1]$ depending only on $d,b,\beta,\Ii$, such that $(X,B)$ has an $(n_0,\Ii_0)$-decomposable $\Rr$-complement. Theorem \ref{thm: betti to complement intro1} follows from Diophantine approximation as in the proof of \cite[Theorem 1.8]{HLS19} (see \cite[Section 6]{HLS19} for details).
\end{proof}

\begin{proof}[Proof of Theorem \ref{thm: betti to complement intro}]
It follows from Theorem \ref{thm: betti to complement intro1}.
\end{proof}

\section{Effective Iitaka fibrations}\label{sec4}
In this section, we prove Theorem \ref{thm: main}.

\begin{lem}\label{lem:horareQ}
Assume that $f: X\rightarrow Z$ is a contraction between normal quasi-projective varieties, and $D$ is an $\Rr$-Cartier $\Rr$-divisor on $X$. Suppose that either
\begin{itemize} 
  \item $\dim Z=0$ and $D\sim_{\Rr}0$, or
  \item $\dim Z>0$ and $D\sim_{\Rr}f^*D_Z$ for some big $\Rr$-Cartier $\Rr$-divisor $D_Z$ on $Z$.
\end{itemize}
Then $\kappa(X,D)\geq 0$ if and only if $D^{h}$ is a $\Qq$-divisor.
\end{lem}
\begin{proof}
We first assume that $\kappa(X,D)\ge0$. Suppose on the contrary that $D^{h}$ is not a $\Qq$-divisor. Let $F$ be a very general fiber of $f$. Then $\{mD|_F\}=\{mD^{h}|_F\}\not=0$ for any positive integer $m$. By our assumption, $D|_F\sim_{\Rr}0$ and
$$\lfloor mD|_F\rfloor=mD|_F-\{mD|_F\}\sim_{\Rr}-\{mD|_F\}$$
is not pseudo-effective for any positive integer $m$. Thus $\lfloor mD\rfloor$ is not pseudo-effective for any positive integer $m$, which implies that $\kappa(X,D)=-\infty$, a contradiction. Therefore, $D^{h}$ is a $\Qq$-divisor.

Now suppose that $D^{h}$ is a $\Qq$-divisor. If $\dim Z=0$, then $D\sim_{\Qq}0$ and thus $\kappa(X,D)=0$. Assume that $\dim Z>0$. Let $b$ be a positive integer, such that $bD\sim 0$ on the generic fiber of $f$. Then there exists $\alpha\in\mathcal{K}(X)$ such that $bD+(\alpha)$ is vertical over $Z$. Since $bD+(\alpha)\sim_{\Rr,Z}0$, $bD+(\alpha)=f^*D_Z'$ for some $\Rr$-Cartier $\Rr$-divisor $D_Z'$ on $Z$ by \cite[Lemma 2.11]{Li20}. Thus $D'_Z\sim_{\Rr}bD_Z$, so $D'_Z$ is big. By Lemma \ref{lem: pullback three iitaka dimensions},
$$\kappa(X,D)=\kappa(X,bD)=\kappa(X,bD+(\alpha))=\kappa(Z,D_Z')=\dim Z>0,$$
and we are done.
\end{proof}

\begin{proof}[Proof of Theorem \ref{thm: main}]
Let $\Ii\subset[0,1]$ be a DCC set. Without loss of generality, we may assume that $1\in\Ii$. Assume that $(X,B)$ is a projective lc pair of dimension $d$ such that $\kappa(X,K_X+B)\geq 0$ and $B\in\Ii$. Since we assume Conjecture \ref{conj: exist gmm} in dimension $d$, $(X,B)$ has a good minimal model $(X',B')$. Possibly replacing $(X,B)$ with $(X',B')$, we may assume that $K_X+B$ is semi-ample. Let $f:X\to Z$ be the ample model of $K_X+B$.

Suppose that $\dim Z=0$. Then $B\in\Ii_0$ for some finite set $\Ii_0\subset\Ii\cap\Qq$ which only depends on $d$ and $\Ii$ by the global ACC \cite[Theorem D]{HMX14} and Lemma \ref{lem:horareQ}, and $K_X+B\sim_\Rr0$. As Conjecture \ref{conj: boundedness and existence of n complement nft} holds in dimension $d$, we may find a positive integer $m_0$ depending only on $d$ and $\Ii_0$ such that $m_0(K_{X}+B)\sim0$. In what follows, we may assume that $\dim Z>0$.

According to Proposition \ref{prop:cbfindex}, there exist a positive integer $p$ and a DCC set $\Ii'\subset[0,1]$ depending only on $d$ and $\Ii$, and a choice of the moduli part $\bM_f$ of the canonical bundle formula for $(X,B)$ over $Z$, such that $B_{Z}\in\Ii'$, $p\bM_f$ is b-Cartier, and
$$p(K_X+B)\sim pf^*\left(K_Z+B_{Z}+\bM_{f,Z}\right),$$
where $B_{Z}$ is the discriminant part of the canonical bundle formula for $(X,B)$ over $Z$. Recall that $Z$ is the ample model of $K_X+B$, hence $K_Z+B_Z+\Mm_{f,Z}$ is big. By \cite[Theorem 1.3]{BZ16}, there is a positive integer $m_1$ depending only on $d,p$, and $\Ii'$, such that $|\lf m_1(K_Z+B_Z+\bM_{f,Z})\rf|$ defines a birational map. By \cite[II 2.11 Lemma]{Nak04}, we have that
\begin{align*}
    H^0\left(X,\lf pm_1(K_X+B)\rf\right)&=H^0\left(X,\lf pm_1f^*\left(K_Z+B_Z+\Mm_{f,Z}\right)\rf\right)\\&\cong H^0(Z,\lf pm_1(K_Z+B_Z+\bM_{f,Z})\rf).
\end{align*}
Therefore $|\lf pm_1(K_X+B)\rf|$ defines a map which is birational to $f_\infty$.

Let $m:=pm_0m_1$, then $m$ satisfies our required property.
\end{proof}

\begin{proof}[Proof of Corollary \ref{cor: ei dim 3}]
It follows from Theorem \ref{thm: main}, the existence of good minimal models in dimension $\leq 3$ \cite{KMM87,KMM94}, and Theorem \ref{thm:3folddcc}.
\end{proof}

\section{Existence of decomposable Iitaka fibrations}\label{sec5}
In this section, we recall the definition of \emph{invariant Iitaka fibrations}, which generalizes Iitaka fibrations to the category of all pairs with non-negative invariant Iitaka dimensions. We then show the existence of decomposable Iitaka fibrations.

\begin{defn}[Invariant Iitaka fibrations]\label{defn: invariant iitaka fibration}
Let $X$ be a normal projective variety and $D$ an $\Rr$-Cartier $\Rr$-divisor on $X$ such that $\kappa_{\iota}(X,D)\geq 0$. A morphism $f_\infty: X_{\infty}\rightarrow Z_{\infty}$ between smooth varieties is called an \emph{invariant Iitaka fibration of $D$} if there exists an $\Rr$-Cartier $\Rr$-divisor $D'$ on $X$, such that $D\sim_{\Rr}D'$, $\kappa(X,D')\geq 0$, and $f_\infty$ is an Iitaka fibration of $D'$.
\end{defn}

By \cite[Lemma 2.3]{Hu20}, an invariant Iitaka fibration of $D$ always exists, and is independent of the choice of $D'$.

We propose the following conjecture, which is a little stronger than Conjecture \ref{conj: decomposable ei intro}.
\begin{conj}\label{conj: decomposable ei strong}
Let $d$ be a positive integer and $\Ii\subset[0,1]$ a DCC set. Then there exist a positive integer $m$, a finite set $\Ii_0\subset (0,1]$, and a DCC set $\Ii'\subset [0,1]$ depending only on $d$ and $\Ii$ satisfying the following. Assume that $(X,B)$ is an lc pair of dimension $d$ such that $B\in\Ii$, $\kappa_\iota(X,K_X+B)\ge0$, and either $\Ii$ is a finite set or all component of $B$ are $\Qq$-Cartier. Then \begin{enumerate}
    \item (Weak version) $(X,B)$ has a $(\Ii_0,\Ii')$-decomposable Iitaka fibration.
    \item (Strong version) $(X,B)$ has an $(m,\Ii_0,\Ii')$-decomposable Iitaka fibration.
\end{enumerate}
\end{conj}

\begin{thm}\label{thm:dif1 strong}
Let $d$ be a positive integer. Assume that the non-vanishing conjecture (Conjecture \ref{conj: non-vanishing}) holds in dimension $d$. Then:
\begin{enumerate}
\item Conjecture \ref{conj: decomposable ei strong}(1) holds in dimension $d$.
\item Assume the effective log Iitaka conjecture (Conjecture \ref{conj: ei intro}) holds in dimension $d$. Then  Conjecture \ref{conj: decomposable ei strong}(2) holds in dimension $d$.
\end{enumerate}
\end{thm}

\begin{proof}[Proof of Theorem \ref{thm:dif1 strong}]
Let $\Ii\subset [0,1]$ be a DCC set and possibly replacing $\Ii$ with $\Ii\cup\{1\}$, we may assume that $1\in\Ii$. 

Assume that $(X,B)$ is an lc pair of dimension $d$ such that $B\in\Ii$, $\kappa_{\iota}(X,K_X+B)\geq 0$, and either $\Ii$ is a finite set or every component of $B$ is $\Qq$-Cartier. Let $f_\infty: X_\infty\to Z_\infty$ be an invariant Iitaka fibration of $K_X+B$, $h: X_\infty\rightarrow X$ the induced morphism, and $F$ a very general fiber of $f_\infty$. Let $B_\infty:=h_*^{-1}B+E$, where $E$ is the sum of all the $h$-exceptional prime divisors. Then $(X_\infty,B_\infty)$ is log smooth, $\kappa_{\iota}(X_\infty,K_{X_\infty}+B_\infty)=\kappa_{\iota}(X,K_X+B)$, and $\kappa_{\iota}(F,(K_{X_\infty}+B_\infty)|_F)=0$.

If $\Ii$ is a finite set, then we let $\tilde\Ii:=\Ii$, otherwise we let $\tilde\Ii:=\emptyset$. We let $\bm{v}_0:=(v^0_1,\dots,v^0_m),g$, and $V$ be as in \cite[Proposition 5.1]{CHL22} which depend only on $d,\Ii,\tilde\Ii$. We may write $B_\infty=\sum b_jB_{\infty}^{(j)}$, where $B_{_\infty}^{(j)}$ are the irreducible components of $B_\infty$. Then there exist distinct Weil divisors $B_{\infty,1},\dots,B_{\infty,m}\ge0$ on $X_\infty$, such that
\begin{itemize}
\item $g(\gamma)\geq\gamma$ for any $\gamma\in\bar\Ii$, and $g(\gamma')=\gamma'$ for any $\gamma'\in\tilde\Ii$,
\item $B_\infty(\bm{v}_0)=\sum g(b_j)B_\infty^{(j)}$, where $B_\infty(\bm{v}):=\sum_{i=1}^m v_iB_{\infty,i}$ for any $\bm{v}:=(v_1,\dots,v_m)\in\Rr^m$,
\item both $(X_\infty,B_\infty(\bm{v}))$ and $(X_\infty,B_\infty(\bm{v})-D_\infty)$ are lc for any $\bm{v}\in V$, where $D_\infty:=B_\infty(\bm{v_0})-B_\infty\ge0$, and
\item $\kappa(X_\infty,K_{X_\infty}+B_\infty(\bm{v})-D_\infty)=\kappa_{\iota}(X_\infty,K_{X_\infty}+B_\infty)=\dim Z$ for any $\bm{v}\in V\cap\Qq^m$.
\end{itemize}
Let $D_F:=D_\infty|_F$, $D:=h_*D_\infty$, $B_F(\bm{v}):=B_\infty(\bm{v})|_F$ and $B(\bm{v}):=h_*B_\infty(\bm{v})$ for any $\bm{v}\in\Rr^m$. By the construction of $\bm{v}_0,g$, and $V$, we have
\begin{itemize}
    \item $(X,B(\bm{v})-D)$ is lc for any $v\in V$, and
    \item $\kappa(F,K_F+B_F(\bm{v})-D_F)=\kappa_{\iota}(F,(K_{X_\infty}+B_\infty)|_F)=0$ for any $\bm{v}\in V\cap\Qq^m$.
\end{itemize}
Note that if $\Ii$ is a finite set, then $D_\infty=0$ and $D_F=0$. 
By the definition of Iitaka fibrations, $f_\infty$ is an Iitaka fibration of $K_{X_\infty}+B_\infty(\bm{v})-D_\infty$ for any $\bm{v}\in V\cap\Qq^m$. As $(X,B(\bm{v})-D)$ is lc for any $\bm{v}\in V$, 
$$K_{X_\infty}+B_\infty(\bm{v})-D_\infty-h^*(K_X+B(\bm{v})-D)\ge0$$
and is $h$-exceptional for any $\bm{v}\in V$. Therefore $f_\infty$ is an Iitaka fibration of $K_X+B(\bm{v})-D$ for any $\bm{v}\in V\cap\Qq^m$. 

Let $\bm{v}_1,\dots,\bm{v}_k\in V\cap\Qq^m$ be rational points, such that $\bm{v}_0$ is contained in the interior of the convex hull of $\bm{v}_1,\dots,\bm{v}_k$. There exist a DCC set $\Ii'\ni 1$ and a finite set $\Ii_0:=\{a_1,\dots,a_k\}\subset (0,1]$ depending only on $d$ and $\Ii$ such that $B_i:=B(\bm{v}_i)-D\in\Ii'$ for any integer $1\le i\le k$, and $\sum_{i=1}^ka_i\bm{v}_i=\bm{v}_0$. In particular, 
$$K_X+B=\sum_{i=1}^ka_i\left(K_X+B_i\right).$$ 
Recall that $f_\infty$ is an Iitaka fibration of $K_X+B(\bm{v})-D$ for any $\bm{v}\in V\cap\Qq^m$. Thus $f_\infty$ is an Iitaka fibration of $K_X+B_i$ for any integer $1\le i\le k$. We conclude that $(X,B)$ has a $(\Ii_0,\Ii')$-decomposable Iitaka fibration, which is (1).

Now suppose that Conjecture \ref{conj: ei intro} holds in dimension $d$. Then there exists a positive integer $m$ depending only on $d$ and $\Ii'$, such that the map defined by $|\lfloor m(K_X+B_i)\rfloor|$ is birational to $f_\infty$ for any integer $1\le i\le k$. Therefore $(2)$ holds.
\end{proof}

\begin{cor}\label{cor: decomposable ei dim 3 strong}
Conjecture \ref{conj: decomposable ei strong} holds when $d\leq 3$.
\end{cor}
\begin{proof}
It follows from Theorem \ref{thm:dif1 strong} and the existence of good minimal models in dimension $\leq 3$ \cite{KMM87,KMM94}.
\end{proof}

\begin{thm}\label{thm:dif2 strong}
Let $d$ be a positive integer. Assume that the good minimal model conjecture (Conjecture \ref{conj: exist gmm}) and the existence of complements (Conjecture \ref{conj: boundedness and existence of n complement nft}) hold in dimension $d$. Then Conjecture \ref{conj: decomposable ei strong} holds in dimension $d$.
\end{thm}
\begin{proof}
It follows from Theorems \ref{thm: main} and \ref{thm:dif1 strong}.
\end{proof}

\begin{proof}[Proof of Theorem \ref{thm:dif1}]
It is a special case of Theorem \ref{thm:dif1 strong}.
\end{proof}

\begin{proof}[Proof of Corollary \ref{cor: decomposable ei dim 3}]
It is a special case of Corollary \ref{cor: decomposable ei dim 3 strong}.
\end{proof}

\begin{proof}[Proof of Theorem \ref{thm:dif2}]
It is a special case of Theorem \ref{thm:dif2 strong}.
\end{proof}

\end{document}